\documentclass[openbib,reqno,12pt]{article}

\advance\textwidth by 1.2in \advance\oddsidemargin by -.6in
\advance\evensidemargin by -.6in
\usepackage{mathrsfs,upgreek}
\usepackage[all,cmtip]{xy}
\usepackage{amsfonts,amssymb,latexsym,amsthm}
\usepackage{amsmath}
\usepackage{hyperref}
\usepackage{amscd, enumerate}
\usepackage{verbatim}

\theoremstyle{definition}
\numberwithin{equation}{section}
\newtheorem{st}{statement}[section]
\newtheorem{lem}[st]{Lemma}

\newtheorem{prop}[st]{Proposition}
\newtheorem{rep.prop}[st]{Proposition}
\newtheorem{thm}[st]{Theorem}
\newtheorem{cor}[st]{Corollary}

\theoremstyle{definition}
\newtheorem{rep.lem}[st]{Lemma}

\newtheorem{defn}[st]{Definition}
\newtheorem{para}[st]{}

\newtheorem{rem}[st]{Remark}
\newtheorem*{notn}{Notation}

\newtheorem*{pf}{Proof}

\newtheorem*{key}{\it{Keywords}}

\newcommand{\eps}{\epsilon}

\def\lar{\mathbb{C}[t^{\pm r}]}

\def\wt{\widetilde} 

\def\fd{finite dimensional\ }
\def\iff{if and only if\ }

\def\lalg{$\mathfrak{g}_{_A} $}
\def\tlalg{${\tilde{\mathfrak{g}}}_{_A} $}

\def\th{$\tilde{\mathfrak{h}} \ $}
\def\hl{$\mathfrak{h}_{_A}$}
\def\thl{${\tilde{\mathfrak{h}}}_{_A}$}

\def\twist{$\mathfrak{g}_{_A}(\mu)$}
\def\ttwist{${\tilde{\mathfrak{g}}}_{_A}(\mu) $}
\def\tV{$V^{(\mu)} \ $}

\def\thlb{{\tilde{\mathfrak{h}}}_{_A}  }

\def\hlb{\mathfrak{h}_{_A} \ }
\def\lalgb{\mathfrak{g}_{_A} \ }
\def\tlalgb{{\tilde{\mathfrak{g}}}_{_A} \ }

\def\nplb{\mathfrak{n}^+_{_A}  }
\def\nnlb{\mathfrak{n}^-_{_A}  }

\def\twistb{\mathfrak{g}_{_A}(\mu) \ }
\def\ttwistb{{\tilde{\mathfrak{g}}}_{_A}(\mu) }
\def\Ifin{$I_{fin}$ }
\def\Ifinmu{$I_{fin}^{\mu}$ }

\def\Bb{\mathbb} 
\def\frk{\mathfrak} 
\def\cal{\mathcal} 
\def\bm{\boldmath}
\def\lb{\lbrace} 
\def\rb{\rbrace}
 
\def\und{\underline}

\def\lamu{\mbox{$_{(\lambda,a,\mu)}$}}
\def\psila{$\psi${\bm{$_{(\lambda,a)}$}}}

\def\psixib{$\psi${\bm{$_{(\xi,b)}$}}}
\def\psixibmu{$\psi${\bm{$_{(\xi,b,\mu)}$}}}
\def\psixibmui{$\psi${\bm{$_{(\xi(i),b,\mu)}$}}}
\def\psixibmus{$\psi${\bm{$_{(\xi(s),b,\mu)}$}}}
\def\psilamu{$\psi{\bm{\lamu}}$}

\def\vpsila{{\bm{$v_{\mbox{\psila}}$}}}

\def\hatm{\hat{m}}

\def\iso{\backsimeq}
\def\A{$A$}

\def\B{Batra, 2004} 
\def\BBFP{Allison, et.all,}
\def\C{Chari, 1986} 
\def\CPone{Chari and Pressley, 1986}
\def\CPtwo{Chari and Pressley, 1988}

\def\K{Kac, 1990} 
\def\R{Rajan, 2004}
\def\Rninethree{Rao, 1993} 
\def\Rninefive{Rao, 1995} 
\def\Rone{Rao, 2001}
\def\Rfour{Rao, 2004}
\def\Y{Yoon, 2002}

\title{Representations of  Graded Multi-Loop \  Lie \ Algebras }
\date{}
\author{ Tanusree Pal and Punita Batra 
\vspace{0.3cm}\\{\small Harish-Chandra Research Institute},\\ {\small Chhatnag 
Road, Jhunsi, Allahabad -211019, India}\\ {\small Email:\ tanusree@hri.res.in,
batra@hri.res.in }}
\begin{document}

\maketitle

\begin{abstract}
Let ``\tlalg (respectively, \ttwist) be the graded multi-loop Lie algebra 
(respectively graded twisted multi-loop Lie algebra)'' associated with the 
simple finite dimensional Lie algebra $\frk{g}$ over $\Bb{C}$. In this paper, 
we prove that irreducible integrable \ttwist-modules with finite dimensional
weight spaces are either highest weight modules or their duals and classify 
the isomorphism classes of irreducible integrable \tlalg-modules and 
\ttwist-modules with finite dimensional weight spaces.\\

\begin{key} Irreducible Integrable module; Multi-loop Lie algebra. \end{key}
\end{abstract}

\section*{\small{INTRODUCTION} }

Let $\mathfrak{g}$ be a finite dimensional simple Lie algebra over 
$\mathbb{C}$. Let A = $\mathbb{C}[t_{1}^{\pm 1} ,\cdots,t_{n}^{\pm 1}]$ be the 
ring of Laurent polynomials in n commuting variables $t_1,t_2,\cdots,t_n$.
Let $\mu_1,\mu_2,\cdots,\mu_n$ be a set of pairwise commuting finite order 
automorphism of $\frk{g}$ such that $\mu_i^{m_i}=Id_{\frk{g}}$ for each i. As 
defined in (\BBFP to appear), a multi-loop Lie algebra $M(\frk{g;\mu_1,\cdots,
\mu_n})$ is given by \begin{equation} M (\frk{g;\mu_1,\cdots,\mu_n}) =
\sum_{(\ell_1,\cdots,\ell_n)\in \Bb{Z}^n} \, 
\frk{g}^{(\bar{\ell_1},\cdots,\bar{\ell_n})} \otimes 
t_1^{\ell_1}\cdots t_n^{\ell_n} \end{equation}
where $\frk{g}^{(\bar{\ell_1},\cdots,\bar{\ell_n})} = \lb x\in \frk{g} \mid 
\mu_ix= \eps_{m_i}^{\ell_i} x \ \rm{for}\ 1\leq i\leq n \rb$. Let $D$ be the 
complex linear span of the derivations $\lb d_1,\cdots,d_n\rb$
of A. A graded multi-loop Lie algebra is defined by
$$M(\frk{g;\mu_1,\cdots,\mu_n}) \oplus D.$$ 
Let \lalg = $M(\frk{g};id_{\frk{g}},\cdots,id_{\frk{g}})$ and \twist =
$M(\frk{g};\mu,id_{\frk{g}},\cdots,id_{\frk{g}})$. 
The universal central extension of \tlalg\ is called a toroidal Lie algebra. 
Integrable irreducible representations for the multi-loop Lie algebras and 
the toroidal Lie algebras have been studied at length in (\B; \Rninefive; 
\Rone; \Rfour). 

The aim of this paper is to classify the isomorphism classes of the irreducible
integrable modules with finite dimensional weight spaces for the graded 
multi-loop Lie algebras
$$ \tlalgb = \lalgb\ \oplus D;  \quad \qquad \ttwistb = \twistb \oplus D.$$
This classification has been given in Theorem 3.8 for \tlalg-modules and 
Theorem 4.7 for \ttwist-modules.
This paper is organized as follows. 

Section 1 reviews known facts about the multi-loop Lie algebras \lalg\ ,\twist 
and the graded multi-loop Lie algebras \tlalg\ ,\ttwist. Let \Ifin and \Ifinmu 
be the categories of integrable modules with finite dimensional weight spaces 
for \tlalg\ and \ttwist\ respectively. It had been proved in (\Rone) that 
irreducible \tlalg\ modules in \Ifin are highest weight modules or their duals.
In Section 2, we show that irreducible \ttwist\ modules in \Ifinmu are either 
highest weight modules or their duals. After recalling the definitions of 
graded and non-graded highest weight modules in the generality of multi-loop 
Lie algebras from (\Rfour), the one-one correspondence between them is used 
in Proposition \ref{restriction} to prove that every irreducible 
highest weight module in \Ifinmu is a \ttwist\ submodule of some irreducible 
highest weight module in \Ifin.

It had been proved in (\Rfour), that an  irreducible highest weight 
\tlalg\ module in \Ifin is of the form ``V(\psila)'', where 
\psila:U($\frk{h}_{_A}$) $\rightarrow A$ is an algebra homomorphism as defined 
in \ref{R-map}. In Section \ref{iso-tlalg}, we analyse the map ``\psila"\ and 
study its effect on V(\psila) in 
Theorem \ref{thm}\footnote{A part of this theorem had 
been proved in (\Y), but the proof here is more elementary.}. Then
using Theorem \ref{thm}, we classify the isomorphism classes of irreducible
modules in \Ifin. In Section \ref{iso-ttwist}, using 
Proposition \ref{restriction} and the general form ``V(\psila)'' of irreducible
integrable highest weight \tlalg\ modules, we give the general form for 
irreducible integrable highest weight \ttwist-modules and classify the 
isomorphism classes of irreducible modules in \Ifinmu. Further it is shown that
an irreducible \tlalg\ module V(\psila)$\in$ \Ifin is completely reducible as a
\ttwist-module if either the corresponding finite dimensional \lalg\ module is
irreducible as a \twist-module or if Image \psila\ = $A$. 

The result on the unique factorization of tensor product of irreducible finite 
dimensional representations of $\frk{g}$, as proved in (\R), considerably 
simplifies the proofs of the main theorems,  Theorem \ref{isomorphism3} and 
Theorem \ref{tiso}, wherein the isomorphism classes of irreducible modules in 
\Ifin and \Ifinmu have been classified.

\noindent{\bf{Notation :}} In the present paper all Lie algebras are defined 
over the complex field $\Bb{C}$.  The superscript $^*$ stands for dual space, 
U($\cdot$) stands for the universal enveloping algebra. $\mathbb{Z}_{+} := 
\lbrace n \in \mathbb{Z} \mid n \geq 0 \rbrace$. For any integer n, $I_n$ 
denotes the set $\lb 1,\cdots,n \rb$. We use underlined letters $\und{m}$ to 
denote a n-tuple of integers ($m_1, m_2, \ldots , m_n$) and in particular for 
every $1\leq i\leq n$, $\und{e}_i$ denotes the n-tuple 
($0,\cdots,\underbrace{1}_{i^{th}},\cdots, 0$) $\in \Bb{Z}^n$ with 1 is in the 
$i^{th}$ place and zero elsewhere. Given an integer p, p$\und{m}$ denotes the 
n-tuple ($pm_1,pm_2,\ldots, pm_n$).

\section {\small{PRELIMINARIES}}
\label{prelims}

\para
Let $\mathfrak{g}$ be a complex finite dimensional simple Lie algebra of rank 
$d$, $\mathfrak{h}$ a Cartan subalgebra of $\mathfrak{g}$ and $\triangle$ 
the set of roots of $\mathfrak{g}$ with respect to $\mathfrak{h}$. Let 
$\pi = \lbrace\ \alpha_1, \alpha_2, \ldots, \alpha_d \ \rbrace$ be a set of 
simple roots. Let $\triangle_+$ (respectively, $\triangle_-$) be the set of 
positive (respectively, negative) roots with respect to $\pi$. 
$\mathfrak{g}$ has a standard decomposition given by:
\[\begin{array}{ll}
\qquad \quad \mathfrak{g} = {\mathfrak{n}}^+ \oplus \mathfrak{h} \oplus 
{\mathfrak{n}}^-, \ \quad \text{where} \quad \  \mathfrak{n}^{\pm} = 
\oplus_{\alpha \in \triangle_{\pm}} \mathfrak{g}_{\alpha}, \quad 
\end{array}\]
and $\mathfrak{g}_{\alpha}$ is the root space corresponding to the root
$\alpha \in \triangle.$ 
For each $\alpha \in \triangle$, there exists $x_{\alpha}^{\pm} \in 
\frk{g}_{\pm \alpha}$ such that [$x_{\alpha}^+,x_{\alpha}^-$]=$\alpha^{\vee}$
and [$\alpha^{\vee},x_{\alpha}^{\pm}$]=$\pm 2 x_{\alpha}^{\pm}$.
 
Let $\mu$ be a $k$-order diagram automorphism of $\frk{g}$. Then 
$\mathfrak{g} = \overset{k-1}{\underset{i=0}{\oplus}} \frk{g}_i$, 
where for $i \in \Bb{Z}$/$k\Bb{Z}$, $\mathfrak{g}_i$ = $\lbrace X \in 
\mathfrak{g} \mid \mu(X) = \epsilon^{i} X \rbrace$, $\epsilon$ is a primitive 
$k^{th}$ root of unity and $\frk{h}_i=\frk{g}_i\cap\frk{h}.$
It has been shown in (\K, Proposition 8.3), that 
the $\mu$ fixed point subalgebra $\mathfrak{g}_0$ of $\mathfrak{g}$ is a 
simple Lie algebra with Cartan subalgebra $\mathfrak{h}_0$ and Chevalley 
generators $x_i^{\pm}, \ i \in$ I, where I is a finite set and $\mathfrak{g}_1$
is an irreducible $\mathfrak{g}_0$ module with lowest weight $\theta_0$
and weight vectors $E_0\in (\frk{g}_1)_{-\theta_0}$ and 
$F_0\in (\frk{g}_1)_{\theta_0}$ .  Let $\triangle_0$ be the root system of 
$\mathfrak{g}_0$.Then $\triangle_0 = (\triangle_0)^s \cup  
(\triangle_0)^{\ell}$, where $(\triangle_0)^{\ell}$ and $(\triangle_0)^{s}$ 
denote the set of long and short roots of $\mathfrak{g}_0$ as given in 
(\K, Chapter 7).

Let \A = $\mathbb{C}[ t_{1}^{\pm 1} ,\ldots, t_{n}^{\pm 1} ]$ be the ring of 
Laurent polynomials in n commuting variables $ t_1,\cdots, t_n$. For 
$\underline{m}$ = ($m_1, m_2, \ldots , m_n$) $\in \mathbb{Z}^n$, we denote
$t_1^{m_1}t_2^{m_2} \ldots  t_n^{m_n}$ by $t^{\underline{m}}$. For a vector 
space V over $\Bb{C}$, let V$_A$ := V$\otimes_{\Bb{C}} A$ and let 
$v(\underline{m}) = v \otimes t^{\underline{m}}$ for $v \in$ V. Via the 
injective map $v\mapsto v\otimes1$ we identify V$\otimes 1$ with the vector 
space V.

Let $\mathfrak{g}_A = \mathfrak{g} \otimes_{_{\Bb{C}}}$ \A  be the multi-loop 
algebra of $\frk{g}$ with Lie bracket given by 
\begin{equation}
\left[ X(\underline{m}) , Y(\underline{n}) \right] = 
\left[ X , Y \right](\underline{m} + \underline{n}), \quad 
\text{for} \ X , Y \in \mathfrak{g},\ \underline{m}, \underline{n} \in 
\mathbb{Z}^n.  \end{equation} Let $D$ be the linear span of the derivations 
$d_1, d_2, \ldots, d_n$ of \A over
$\mathbb{C}$ defined by $d_i(t^{\und{m}})=m_i t^{\und{m}}$ for $1\leq i\leq n$ 
and $\und{m}\in \Bb{Z}^n$. Denote by \tlalg = \lalg\ $\oplus D$ the graded 
multi-loop Lie algebra of $\frk{g}$, in which [$d_i ,  X(\underline{m})$]=
$ m_i X(\und{m})$ for all $\und{m} \in \Bb{Z}^n$ and  $1\leq i \leq n.$ 
\tlalg\ has a decomposition given by :
\[\tlalgb = \nplb \oplus \thlb \oplus \nnlb, \qquad \text{where} \ \thlb = \hlb
\oplus D. \]  

Given a diagram automorphism $\mu$ of a simple Lie algebra $\frk{g}$, it can be
extended to an automorphism of \tlalg\ as follows : \[\begin{array}{rl}
\mu(X(\und{m})) = \eps^{- m_1} \mu(X)(\und{m}) , \quad & \quad 
\mu(d_i) = d_i, \qquad \forall\ i = 1,\cdots, n. \end{array} \]  
The fixed point subalgebra $ \twistb = \lb X \in \lalgb \vert \mu(X) = X \rb$ 
of \lalg\ is called the twisted multi-loop Lie algebra of $\frk{g}$ (see 
\B) and \ttwist = \twist\ $\oplus D$ is the graded twisted multi-loop 
algebra. For a subalgebra $\mathfrak{a}$ of $\mathfrak{g}$, let 
$\mathfrak{a}$($\mu$) $ = \lb X \in \mathfrak{a} \ \vert \ \mu(X) = X \rb$. 
Then \ttwist\ has a decomposition given by :
\[ \ttwistb = \nplb(\mu) \oplus \thlb(\mu) \oplus \nnlb(\mu), 
\qquad \text{where} \ \thlb(\mu) = \frk{h}_{_A}(\mu) \oplus D. \]
Define elements $\delta_i \in (\mathfrak{h} + D)^*$ by $\delta_i {\vert}_{
\mathfrak{h}}$ = 0, $\delta_i(d_j)$ = $\delta_{ij} $ for i,j $\in \lb 1,2,
\ldots,n \rb$. The space ${\mathfrak{h}}^{*}$ (resp. $\mathfrak{h}_0^*$ ) is 
identified with a subspace of (${\mathfrak{h}}+D$)$^{*}$ (resp. 
(${\mathfrak{h}}_0+D$)$^{*}$) by setting $\lambda(d_i)$ = 0 for 
$\lambda \in {\mathfrak{h}}^{*}$\ (resp. $\mathfrak{h}_0^*$).
For $\und{m} = (m_1,m_2, \ldots, m_n) \in \mathbb{Z}^n$, define 
$\delta_{\underline{m}} = \sum_{i=1}^n \  m_i \delta_i$.

\para Let $\triangle$ be the set of roots for 
(\tlalg, $\frk{h}$ ) with the standard base $\pi$ and $\triangle^{\mu}$ be the 
set of roots for (\ttwist, $\frk{h}_0$) with the standard base $\pi^{\mu}$.
Let W(\tlalg) and W(\ttwist) denote the Weyl groups of \tlalg\ and \ttwist\ 
respectively. $\gamma \in\triangle$ (respectively $\gamma \in
\triangle^{\mu}$) is called a real root if $\gamma$ is W(\tlalg)-conjugate to
$\pi$ (respectively, W(\ttwist)-conjugate to $\pi^{\mu}$). We shall denote the 
set of real roots of $\triangle$ (respectively $\triangle^{\mu}$) by 
$\triangle_{re}$ (respectively $\triangle^{\mu}_{re}$) and  
$\triangle_{im}$ = $\triangle\smallsetminus\triangle_{re} $ 
(respectively $\triangle_{im}^{\mu}$ = $\triangle^{\mu}\smallsetminus
\triangle_{re}^{\mu}$) shall denote the set of imaginary roots of \tlalg 
(respectively \ttwist).

The root space corresponding to a root $\alpha$ of \tlalg (resp. of \ttwist)
is denoted by (\tlalg)$_{\alpha}$ (resp. (\ttwist)$_{\alpha}$). For 
$\alpha \in \triangle_{re}$, dim$_{\Bb{C}}$ (\tlalg)$_{\alpha}$=1
(resp. for $\alpha \in \triangle_{re}^{\mu}$, 
dim$_{\Bb{C}}$ (\ttwist)$_{\alpha}$=1 ).

\para
From the theory of \fd \ simple Lie algebras over $\Bb{C}$ we know that given
a root $\alpha \in \triangle$, there exists $x_{\alpha}^{\pm} \in 
\mathfrak{g}_{\pm \alpha}$ and ${\alpha}^{\vee} \in \frak{h}$ such that 
S$_{\alpha}$ = $\Bb{C}$-span $\lb x_{\alpha}^{\pm},{\alpha}^{\vee} \rb$ is 
isomorphic to $sl_2$. It has been shown in (\Rfour), that for such a choice 
of elements $\lb x_{\alpha}^{\pm}, {\alpha}^{\vee} \rb$ in $\mathfrak{g}$, the 
$\Bb{C}$-span of $\lb x_{\alpha}^{+}(\und{m}),x_{\alpha}^{-}(-\und{m}), 
{\alpha}^{\vee} \rb$ forms a copy of $sl_2$ in \tlalg\ for all $\und{m} \in 
\Bb{Z}^n$.

\section{\small{INTEGRABLE IRREDUCIBLE REPRESENTATIONS OF \tlalg\ AND 
\ttwist\\ WITH FINITE DIMENSIONAL WEIGHT SPACES}}
\label{integrable}

\begin{defn} A \tlalg-module V (resp. \ttwist-module \tV) is said to be 
integrable if \begin{verse}
(1). V is a weight module with respect to  \th = $\mathfrak{h} \oplus D$ 
(resp. \tV is a weight module with respect to $\tilde{\frk{h}_0}$ =
$\mathfrak{h}_0 \oplus D$ ), 
\[  V = \bigoplus_{\lambda \in ({\tilde{\mathfrak{h}}})^{*}} 
V_{\lambda} \qquad (\rm{resp.}  V^{(\mu)} = \bigoplus_{\lambda \in 
(\tilde{\frk{h}_0})^*} V_{\lambda}^{(\mu)} ).\]\\
(2). For every $v \in$ V and $\alpha$ $\in \triangle_{re}$ 
(resp. $\alpha \in \triangle^{\mu}_{re}$ ),   $\exists$ a positive integer 
N = N($\alpha,v$) such that $X_{\alpha}^N.v=0$ for all $X_{\alpha} \in 
(\tlalgb)_{\alpha}$ (resp. $X_{\alpha} \in (\ttwistb)_{\alpha}$).
\end{verse}
\end{defn}
Let \Ifin (resp. \Ifinmu ) be the category of of integrable \tlalg-modules 
(resp. \ttwist-modules) with \fd weight spaces. 

\begin{prop} : \label{H.wt}
If V is an integrable irreducible module with finite dimensional
weight spaces for the Lie algebra \tlalg or \ttwist, then there exists a 
non-zero weight vector $v \in$ V such that
\[\begin{array}{r}
{\mathfrak{n}}_{A}^+ v =0, \quad ({\rm{resp}}.  {\mathfrak{n}}_{A}^+(\mu) v=0),
\\ or \quad {\mathfrak{n}}_{A}^- v =0. \quad ({\rm{resp}}. \  
{\mathfrak{n}}_{A}^-(\mu) v =0). \end{array}\] \end{prop}
\pf The proposition was proved for irreducible \tlalg-modules 
V $\in$ \Ifin in (\Rone, Proposition 3.2).
If V $\in$ \Ifinmu is irreducible, then the proof follows from 
(\CPtwo, Proposition 3.5).  $\hfill \blacksquare$ \endproof

\begin{defn}(\Rfour).
A \tlalg-module $V$, is called graded (resp. non-graded) highest weight module 
for \tlalg (resp. \lalg) if there exist a weight vector $v \in $V with respect 
to $\frk{h} \oplus  D$(resp. $\frk{h}$) such that 
\begin{verse}
(1). $V = U (\tlalgb).v$, $\left(\mbox{resp. V=U(\lalg).$v$}\right)$.\\
(2). $\nplb .v=0$.  \\
(3). $U(\thlb).v$ is an irreducible module for \thl, 
$\left(resp. \mbox{ $h v= \Psi (h) v$  $\ \forall \ h \in \hlb$, for  
$\Psi \in (\frk{h}_A)^*$}\right)$.
\end{verse}
\end{defn}

\para
Let $\psi:U(\frk{h}_A$)$ \rightarrow A$  be a $\Bb{Z}^n$-graded algebra
homomorphism. Then $A$ is a module for $U(\frk{h}_A$)  via $\psi$ defined as:
$$\begin{array}{ll}
h(\und{m})t^{\und{s}} = \psi(h(\und{m}))t^{\und{s}}, & h\in \frk{h}, 
\und{m}\in\Bb{Z}^n\\ ht^{\und{s}} = \psi(h)t^{\und{s}}, & h\in\frk{h}.
\end{array}$$ Let $A_{\psi}$ = Image of $\psi$. If $A_{\psi}$ is an 
irreducible $\tilde{\frk{h}}_A$-module, let $\frk{n}_{A}^+$ act trivially on 
$A_{\psi}$. Consider the following induced \tlalg-module 
$M(\psi) = U(\tlalgb) \otimes_{\cal{\tilde{B}}} A_{\psi}$, where 
$\cal{\tilde{B}}= \nplb \oplus \thlb$. As shown in (\Rone, Proposition 1.4),
$M(\psi)$ has a unique irreducible quotient $V(\psi)$.

\noindent The following criteria for irreducibility of $A_{\psi}$ as a 
$\tilde{\frk{h}}_A$-module had been proved in \Rninefive:
\begin{lem}(\Rninefive, Lemma 1.2) \label{Apsi} A$_{\psi}$ is an irreducible 
${\tilde{\mathfrak{h}}}_A$-module if and only if each non-zero
homogeneous element of A$_{\psi}$ is invertible. 
\end{lem}

Let $\Psi\in (\frk{h}_A)^*$ and let $\frk{h}_A$ act on the one dimensional 
vector space $\Bb{C}(\Psi)$ by $\Psi$. Let $\frk{n}_{A}^+$ act trivially on 
$\Bb{C}(\Psi)$. Consider the following induced \lalg-module 
$M(\Psi) = U(\lalgb) \otimes_{\cal{B}} \Bb{C}(\Psi)$, where 
$\cal{B}= \nplb \oplus \hlb$. Then $M(\Psi)$ has a unique irreducible quotient 
$V(\Psi)$. 

With $\psi$ as above, let $\Psi = \cal{E}\cdot\psi:A_{\psi}\rightarrow \Bb{C}
$, where $\cal{E}:A \rightarrow \Bb{C}$ is the evaluation map defined by 
$\cal{E}(t^{\und{m}}) = 1$. $V(\Psi) = V(\cal{E}\cdot\psi)$ can be made into a 
(graded) \tlalg-module by 
$$\begin{array}{ll}
g(\und{m}).v(\und{r}) = (g(\und{m}.v))(\und{m}+\und{r}), & \text{for}\ 
g\in\frk{g}, \und{m},\und{r}\in\Bb{Z}^n, v\in V(\Psi),\\
d_i(v(\und{r})) = r_i v(\und{r}), & \text{for}\ 
\und{r}\in\Bb{Z}^n, v\in V(\Psi), 1\leq i\leq n,\\
hv(\und{r}) = (hv)(\und{r}), & \text{for}\ h\in\frk{h},
\und{r}\in\Bb{Z}^n, v\in V(\Psi). \end{array}$$

\begin{prop}[Proposition 3.5,\Rfour]\label{irreducible comp}  
Let $\psi$ and $\cal{E}\cdot\psi$ be as above. Assume that 
$A_{\psi}$ is an irreducible \thl module. Let $G$ be the set of 
coset representatives for $A/A_{\psi}$. If $v$ be a highest weight vector of 
$V$($\cal{E}\cdot\psi$), then: \\
(1). ${V}(\cal{E}\cdot\psi)\otimes A = \underset{\wt{\und{m}}\in G }{\oplus}\ 
U(\tlalgb)v(\wt{\und{m}})$ as \tlalg modules.  \\
(2). For each $\wt{\und{m}}\in G$, $U(\tlalgb)v(\wt{\und{m}})$ is an 
irreducible \tlalg module.\\
(3). U(\tlalg)$v(\wt{\und{0}}) \iso V (\psi)$ as \tlalg module.
\end{prop}

\begin{lem} \label{fin-d} 
Let V $\in$ \Ifin be an irreducible highest weight \lalg-module such that
$$h.v =\Psi(h)v, \qquad \forall\ h\in \frk{h}_A,$$ where $v$ is the highest 
weight vector of V and $\Psi\in (\frk{h}_A)^*$. Then V is isomorphic to a 
finite dimensional $\frk{g}$ module and $\Psi\vert_{\frk{h}}$ is a dominant
integral weight of $\frk{g}$. \end{lem} 
\pf Follows from (Rao,2004, Lemma 3.6 and Proposition 3.20). 
$\hfill \blacksquare$ \endproof

\noindent We now discuss the corresponding notions for the \twist\ and 
\ttwist-modules. 
\begin{defn} 
A \ttwist-module $W$, is called graded (resp. non-graded) highest weight module
for \ttwist (resp. \twist) if there exist a weight vector $w\in W $ with 
respect to $\frk{h}_0 \oplus  D$(resp. $\frk{h}_0$) such that 
\begin{verse}
(1). $W = U (\ttwistb).w$, $\left(\mbox{resp. W=U(\twist).$w$}\right)$.\\
(2). $\nplb(\mu) .w=0$.  \\
(3). $U(\tilde{\frk{h}}_A(\mu)).w$ is an irreducible module for 
$\tilde{\frk{h}}_A$, (resp. $h w= \Psi^{\mu} (h) w$,  $\ \forall \ h \in 
\frk{h}_A(\mu)$, for $\Psi^{\mu} \in \big(\frk{h}_A(\mu)\big)^*$ ).
\end{verse}
\end{defn}

Let $\psi^\mu:U(\frk{h}_A(\mu)$)$ \rightarrow A$  be a $\Bb{Z}^n$-graded 
homomorphism. Then $A$ is a module for $U(\frk{h}_A(\mu)$)  via $\psi^{\mu}$ 
defined as:
$$\begin{array}{ll}
h(\und{m})t^{\und{s}} = \psi^{\mu}(h(\und{m}))t^{\und{s}}, & h(\und{m})\in 
\frk{h}_A(\mu),\\ ht^{\und{s}} = \psi(h)t^{\und{s}}, & h\in\frk{h}_0.
\end{array}$$ If Image $\psi^{\mu}$ = $A_{\psi^\mu}$ is an 
irreducible $\tilde{\frk{h}}_A(\mu)$-module, let $\frk{n}_{A}^+(\mu)$ act 
trivially on $A_{\psi^{\mu}}$. Consider the following induced \ttwist-module 
$M(\psi^{\mu}) = U(\ttwistb) \otimes_{\cal{\tilde{B^{\mu}}}} A_{\psi^{\mu}}$, 
where $\cal{\tilde{B^{\mu}}}= \nplb(\mu) \oplus \thlb(\mu)$. It is easy to see
that $M(\psi^{\mu})$ has a unique irreducible quotient $V(\psi^{\mu})$.

Let $\Psi^{\mu}\in (\frk{h}_A(\mu))^*$ and let $\frk{h}_A(\mu)$ act on the one 
dimensional vector space $\Bb{C}(\Psi^{\mu})$ by $\Psi^{\mu}$. Let 
$\frk{n}_{A}^+(\mu)$ act trivially on $\Bb{C}(\Psi^{\mu})$. Consider the 
following induced \twist-module $M(\Psi^{\mu}) = U(\twistb) 
\otimes_{\cal{B^{\mu}}} \Bb{C}(\Psi^{\mu})$, where 
$\cal{B^{\mu}}= \nplb(\mu) \oplus \hlb(\mu)$. Then $M(\Psi^{\mu})$ has a unique
irreducible quotient $V(\Psi^{\mu})$ (see \B, Proposition 2.1). 

With $\psi^{\mu}$ as above, let $\Psi^{\mu} = \cal{E}\cdot \psi^{\mu}:
A_{\psi^{\mu}}\rightarrow \Bb{C} $, where $\cal{E}$ is the evaluation map as 
defined above. Then  $V(\Psi^{\mu}) = V(\cal{E}\cdot \psi^{\mu})$ can be made 
into a (graded) \ttwist-module by 
$$\begin{array}{ll}
g(\und{m}).v(\und{r}) = (g(\und{m}.v))(\und{m}+\und{r}), & \text{for}\ 
g(\und{m})\in \twistb, \und{m},\und{r}\in\Bb{Z}^n, v\in V(\Psi^{\mu}),\\
d_i(v(\und{r})) = r_i v(\und{r}), & \text{for}\ 
\und{r}\in\Bb{Z}^n, v\in V(\Psi^{\mu}), 1\leq i\leq n,\\
hv(\und{r}) = (hv)(\und{r}), & \text{for}\ h\in\frk{h}_0,
\ \und{r}\in\Bb{Z}^n,\  v\in V(\Psi^{\mu}). \end{array}$$

\begin{prop}\label{irreducible comp'}  
Let $\psi^\mu$ and $\cal{E}\cdot\psi^\mu$ be as above. Assume that 
$A_{\psi^{\mu}}$ = Image $\psi^\mu$ is an irreducible 
$\tilde{\frak{h}}_A(\mu)$-module. Let $G^{\mu}$ be the set of 
coset representatives for $A/A_{\psi^{\mu}}$. Let $w$ be a highest weight 
vector of $V$($\cal{E}\cdot\psi^\mu$). Then \\
(1). V$(\cal{E}\cdot\psi^{\mu})\otimes A$ = $\underset{\wt{\und{m}}\in G^{\mu}}
{\oplus}$ U(\ttwist).$w(\und{{m}})$ at \ttwist-modules.\\ 
(2). For each $\wt{\und{m}}\in G^{\mu}$, U(\ttwist)$w(\wt{\und{m}})$ is an 
irreducible \ttwist-module.\\
(3). U(\ttwist).$w(\und{0})$ isomorphic to V($\psi^{\mu}$) as a \ttwist-module.
\end{prop}
\pf Follows on similar lines to the proof of (\Rninethree, Proposition 4.8)
and (\Rninefive, Proposition 1.8).
$\hfill \blacksquare$ \endproof

\begin{lem}\label{twist1} Let \tV $\in$ \Ifinmu be an irreducible \twist-module
such that \begin{eqnarray}
\frk{n}_{_A}^+(\mu).w=0,  \label{1}  \qquad \qquad  \qquad \qquad 
\\
h.w = \Psi^{\mu}(h)w, \qquad \forall\ h \in \frk{h}_{_A}(\mu),\label{2}
\end{eqnarray}where $w$ is the highest 
weight vector of \tV and $\Psi^{\mu} \in  (\frk{h}_{_A}(\mu))^*$.
Then \tV is isomorphic to a finite dimensional $\frk{g}$-module and 
$\Psi^{\mu}\vert_{\frk{h}_0}$ is a dominant integral weight of $\frk{g}$.
\end{lem}
\pf \tV being an irreducible \twist-module with finite-dimensional weight 
spaces, by (Batra, 2004, Theorem 3.5), \tV is a module for the 
finite-dimensional Lie algebra \twist/(\twist I) := \twist(I), where I is a 
cofinite ideal in $A_0 =\Bb{C}[t_1^{\pm k},t_2^{\pm},\cdots,t_n^{\pm}]$.  
However under the conditions Eq (\ref{1}) and Eq (\ref{2}), \tV is a non-graded
highest weight module. Consequently \tV = $U(\nnlb(\mu))(I).w$, where $w$ is 
the highest weight vector of \tV. Let $y_1,y_2,\cdots, y_N$ be a vector basis 
for the finite dimensional Lie algebra \twist(I). By definition of 
integrability, $y_1,y_2,\cdots, y_N$ act as locally finite endomorphisms on 
\tV. Hence using PBW theorem, we conclude \tV is finite dimensional. By 
(Batra, 2004, Theorem 3.7), a finite dimensional \twist-module, is 
isomorphic as a \twist\  module to the $\frk{g}$-module 
$\otimes_{i=}^N \ V(\lambda_i)$, for $\und{\lambda}$ = 
$(\lambda_1,\lambda_2,\cdots,\lambda_N)$ dominant integral weights,
where for each $i$, $V(\lambda_i)$ denotes the irreducible $\frk{g}$- module 
with highest weight $\lambda_i$. Since \tV is finite dimensional,
\tV is isomorphic to $\otimes_{i=1}^N \ V(\lambda_i)$ and hence the lemma.
$\hfill \blacksquare$ \endproof 

By Lemma \ref{fin-d} and Lemma \ref{twist1}, non-graded highest weight 
\lalg-modules and \twist-modules are highest weight $\frk{g}$-modules of the 
form $V(\lambda_1)\otimes\cdots\otimes V(\lambda_N)$, for $\lambda_i$'s 
dominant integral and hence they are finite dimensional irreducible 
$\underset{N-copies}{\oplus}\frk{g}$ modules for some positive integer $N$.

\begin{lem}\label{rank} Let $\psi$ :U(\hl)$\rightarrow A$ be a $\Bb{Z}^n$ 
graded algebra homomorphism such that Image $\psi$ is an irreducible \hl\ 
module and V($\cal{E}\cdot\psi$) is a non-trivial irreducible highest weight
\lalg-module. Let $\Gamma = \lb \und{m}\in \Bb{Z}^n \mid t^{\und{m}} \in 
A_{\psi}\rb$, where $A_{\psi}$ = Image $\psi$. Then $\Gamma$ is a subgroup of 
$\Bb{Z}^n$ of rank n.\end{lem}
\pf Given  $A_{\psi}$ is an irreducible \hl\ module. Since $\psi$ is an algebra
homomorphism, $\Gamma$ is closed under addition. Clearly, for $h\in \frk{h}$,
$\psi(h\otimes 1) = \cal{E}\cdot\psi(h\otimes1)$, where $\cal{E}$ is the 
evaluation map as defined above. But if V($\cal{E}\cdot\psi$) is a non-trivial 
irreducible highest weight \lalg-module then by Lemma \ref{fin-d}, 
$(\cal{E}\cdot\psi)\vert_{\frk{h}}$ is dominant integral and hence there exists
$h\in\frk{h}$ such that $\psi(h\otimes1)\neq 0$. Consequently, 
$\und{0}\in \Gamma$. By Lemma \ref{Apsi}, $A_{\psi}$ every homogeneous element 
of $A_{\psi}$ is invertible. Hence $\Gamma$ is a subgroup $\Bb{Z}^n$.

\noindent Claim: Rank $\Gamma$ = n.\\
On the contrary suppose that rank $\Gamma$ = k$< n$. Then there exists a basis 
vector $\und{a}$ of $\Bb{Z}^n$ such that  
\begin{eqnarray}\psi(\frk{h}\otimes\Bb{C}t^{\und{a}}[t^{\und{a}}])=0, \quad & 
\quad\psi(\frk{h}\otimes\Bb{C}t^{-\und{a}}[t^{-\und{a}}])=0.\end{eqnarray} 
Without loss of generality we may assume $\und{a}=\und{e}_n$. Then it follows 
that in the non-graded \lalg-module V($\cal{E}\cdot\psi$)
\begin{eqnarray}\cal{E}\cdot\psi(\frk{h}\otimes P(t_n,t_n^{-1}))=0, \quad &
\quad \text{for any polynomial $P(t_n,t_n^{-1})\in\Bb{C}[t_n,t_n^{-1}]$}. 
\end{eqnarray}
By Lemma \ref{fin-d}, V($\cal{E}\cdot\psi$) is a finite dimensional 
\lalg-module and there exists Lie algebra homomorphism 
\[\Phi_{fin}:\lalgb\rightarrow {\rm{End}}(V(\cal{E}\cdot\psi)),\] such that 
$\Phi_{fin}(\lalgb)$ is isomorphic to $\underset{N-copies}{\oplus} \frk{g}$ for
some integer $N$.

\noindent Claim : $\frk{h}\otimes t_n$, $\frk{h}\otimes t_n^{-1} \subset$ 
Ker $\Phi_{fin}$.\\
Let V($\cal{E}\cdot\psi$) be the $\frk{g}$-module generated by the highest 
weight vector $v$. Any general element of V($\cal{E}\cdot\psi$) is of the
form $x_{\alpha_1}^-(\und{m}_1)\cdots x_{\alpha_k}(\und{m}_k).v$. We prove the 
claim by showing that $h(\und{e}_n).(x_{\alpha_1}^-(\und{m}_1)\cdots 
x_{\alpha_k}^-(\und{m}_k).v)=0$ for  $\lb \alpha_i\rb_{1\leq i\leq k}\in 
\triangle^+$. This is proved by inducting on $k$.

\noindent For $k=1$, \ \  $h(\und{e}_n)x_{\alpha_1}^-(\und{m}_1)v =
-\alpha_1(h)x_{\alpha_1}^-(\und{e}_n+\und{m}_1)v,$ for $h\in \frk{h}.$\\
Thus \ \ \ \ \ \ \ \  \ $h(\und{e}_n)x_{\alpha_1}^-(\und{m}_1)v =0 \quad 
\text{if and only if} \ \quad x_{\alpha_1}^-(\und{e}_n+\und{m}_1)v = 0.$\\
Observe that, $$x_{\alpha_j}^+x_{\alpha_1}^-(\und{e}_n+\und{m}_1).v =
x_{\alpha_1}^-(\und{e}_n+\und{m}_1)x_{\alpha_j}^+v + 
\delta_{j,1} \alpha_1^{\vee}(\und{e}_n+\und{m}_1)v =0,$$
by assumption and the fact that $v$ is a highest weight vector in 
V($\cal{E}\cdot\psi$). This implies that for a positive root $\alpha_1$ and
$\und{m}_1\in\Bb{Z}^n$, $x_{\alpha_1}^-(\und{e}_n+\und{m}_1).v$ is a highest 
weight vector of V($\cal{E}\cdot\psi$) of weight less than 
$\cal{E}\cdot\psi\vert_{\frk{h}}$, which is a contradiction. Hence 
$x_{\alpha_1}^-(\und{e}_n+\und{m}_1).v = 0$. Thus the claim holds for $k=1$.

Now suppose the claim holds for all $k\leq s-1$.\\
For $k=s$,  $h(\und{e}_n).(x_{\alpha_1}^-(\und{m}_1)\cdots 
x_{\alpha_s}^-(\und{m}_s).v) = x_{\alpha_1}^-(\und{m}_1).h(\und{e}_n)\cdots 
x_{\alpha_k}^-(\und{m}_s).v -\alpha_1(h) x_{\alpha_1}^-(\und{m}_1+\und{e}_n)
\cdots x_{\alpha_s}^-(\und{m}_s).v.$ By induction hypothesis the first term
in the RHS is zero, hence the claim holds for $k=s$ if and only if
$x_{\alpha_1}^-(\und{m}_1+\und{e}_n)\cdots x_{\alpha_s}^-(\und{m}_s).v.=0$.

Applying $x_{\alpha_j}^+$ to 
$x_{\alpha_1}^-(\und{m}_1+\und{e}_n)\cdots x_{\alpha_s}^-(\und{m}_s).v.$ and
using induction hypothesis recursively we get,
$$\begin{array}{l}
x_{\alpha_j}^+x_{\alpha_1}^-(\und{m}_1+\und{e}_n)\cdots 
x_{\alpha_s}^-(\und{m}_s).v.\\ = \delta_{1,j} 
\alpha_1^{\vee}(\und{m}_1+\und{e}_n)x_{\alpha_2}^-(\und{m}_2)\cdots 
x_{\alpha_k}^-(\und{m}_s).v.
+ x_{\alpha_1}^-(\und{m}_1+\und{e}_n) 
x_{\alpha_j}^+\cdots x_{\alpha_s}^-(\und{m}_s).v.\\
= x_{\alpha_1}^-(\und{m}_1+\und{e}_n) 
x_{\alpha_j}^+\cdots x_{\alpha_s}^-(\und{m}_s).v.\\
\vdots\\
=x_{\alpha_1}^-(\und{m}_1+\und{e}_n)\cdots x_{\alpha_j}^+x_{\alpha_s}^-
(\und{m}_s).v.
= x_{\alpha_1}^-(\und{m}_1+\und{e}_n)\cdots x_{\alpha_s}^-
(\und{m}_s) x_{\alpha_j}^+ v=0,
\end{array}$$ since $v$ is a highest weight vector of V($\cal{E}\cdot\psi$).
Now the same weight argument as in the case $k=1$ shows that
$x_{\alpha_1}^-(\und{m}_1+\und{e}_n)\cdots x_{\alpha_s}^-(\und{m}_s).v.=0$.
Hence the claim holds for $k=s$. 

Thus we get that $\frk{h}\otimes t_n$, $\frk{h}\otimes t_n^{-1} \subset$ 
Ker $\Phi_{fin}$. But Ker $\Phi_{fin}$ is of the form $\frk{g}\otimes I$ for 
some ideal $I$ of $A$ and the smallest ideal of $A$ containing $t_n$ and 
$t_n^{-1}$ is $A$ itself. Thus rank $\Gamma < n$, would imply that 
V($\cal{E}\cdot\psi$) is a trivial \lalg-module. Hence the lemma. 
$\hfill \blacksquare$ \endpf

\begin{prop}\label{restriction}
Every irreducible \ttwist-module \tV $\in$ \Ifinmu is a \ttwist-submodule of 
some irreducible \tlalg-module V $\in$ \Ifin. \end{prop}
\pf  Let V($\psi^{\mu}$) be an irreducible \ttwist-module, where 
$\psi^{\mu}:$ U($\frk{h}_{_A}(\mu)$) $\rightarrow A$ is an algebra 
homomorphism. Then by Lemma \ref{twist1}, V$(\cal{E}\cdot\psi^{\mu})$ is a 
finite dimensional \twist-module of the form $V$({\bm{$\lambda$}}):= 
$V(\lambda_1)\otimes\cdots\otimes V(\lambda_N)$ for some positive integer N. 
The action of \twist\ on $V$({\bm{$\lambda$}}) is given by a surjective 
homomorphism
$$\Phi_{fin}^{\mu}:\frk{g}_{_A}(\mu)\rightarrow \underset{N-copies}{\oplus} 
\frk{g}.$$ From Lemma \ref{twist1}, it is clear that $\Phi_{fin}^{\mu}$ can be 
extended to give a surjective map 
$$\Phi_{fin}:\frk{g}_{_A}\rightarrow \underset{N-copies}{\oplus} \frk{g},$$
and thus give $V$({\bm{$\lambda$}}) the structure of an irreducible non-graded 
\lalg-module. 

Let $w$ be the highest weight vector of $V$({\bm{$\lambda$}}). Then by 
Proposition \ref{irreducible comp} and Proposition \ref{irreducible comp'}, 
U(\tlalg)$w(\und{0})$ is an irreducible \tlalg\ module and the irreducible 
\ttwist\ module V($\psi^{\mu}$), is the \ttwist\ irreducible component 
U(\ttwist).$w(\und{0})$ of $V$({\bm{$\lambda$}})$_{A}$. But 
U(\ttwist).$w(\und{0}) \subset$ U(\tlalg).$w(\und{0})$. Hence the proposition. 
$\hfill \blacksquare$ \endproof

\section{\small{ISOMORPHISM CLASSES OF IRREDUCIBLE \tlalg-MODULES IN \Ifin}}
\label{iso-tlalg}
It had  been proved in (\Rfour) :
\begin{prop}(\Rfour, Proposition 3.20) \label{R-map} Suppose $V(\psi)$ $\in$
\Ifin is an  irreducible highest weight module for \tlalg.  Then the 
$\Bb{Z}^n$ graded algebra homomorphism 
$\psi$= \psila:U($\frk{h}_{_A}$) $\rightarrow A$ is given by
\[\mbox{\psila} (h \otimes t^{\underline{m}}) =
\displaystyle{\sum_{j=0}^{N}} a_{I_j}^{\underline{m}} \lambda_j (h)
t^{\underline{m}}, \ {\rm for} \ h \in \frk{h},\]
where for each $1\leq i\leq n$, ${\bf{a}}_i = (a_{i1}, \cdots, a_{i N_i})$
and $a_{ij}$'s are distinct non-zero complex numbers for $1\leq j\leq N_i$, 
$a_{I_j}^{\underline{m}}$ are as defined in {\bf{\ref{construct}}} and each 
$\lambda_j$ is a dominant integral weight of $\frk{g}$.
\end{prop}

In this section we shall first recall the definition of the $\Bb{Z}^n$-graded
algebra homomorphism  \psila : $U(\hlb)\rightarrow A$ as given in 
(\Rfour, Section 3) and then give a the isomorphism classes of the graded 
irreducible integrable \tlalg-modules V(\psila).

\para\label{construct}{\it
{The construction of the $\Bb{Z}^n$ graded algebra homomorphism 
$\psi: U(\hlb)\rightarrow A$.}} 
Let $n$ be a positive integer.  For each $i,  1 \leq i \leq
n$, let $N_i$ be a positive integer.  Let $\underline{a}_i = 
(a_{i1}, \cdots, a_{i N_i})$ be non-zero distinct complex numbers.  
Let $N=N_1 \cdots N_n$.  Let $I = (i_1, \cdots, i_n)$ where 
$1 \leq i_j \leq N_j$.  For $\underline{m} = (m_1, \cdots, m_n) \in \Bb{Z}^n$, 
define $a_I^{\underline{m}} = a_{1i_1}^{m_1} \cdots a_{ni_n}^{m_n}$. For a pair
({\bm{$\lambda,a$}}) $\in ({\mathfrak{h}}^{*})^N \times 
(\Bb{C}^{\times})^{N}$, define $\psi${\bm{$_{(\lambda,a)}$}} :U($\hlb$) 
$\rightarrow$ A  by extending 
\begin{eqnarray}\label{eqndef}
\psi{\mbox{\bm{$_{(\lambda,a)}$}}}
(h \otimes t^{\underline{m}}) = \left(\sum_{i=1}^N \  {\lambda}_I(h) \ 
a_{I}^{\und{m}} \right)  t^{\und{m}}, \  & \forall \ h \in \frk{h}, \und{m} \in
\Bb{Z}^n, \end{eqnarray}
to an algebra homomorphism. Clearly \psila\ is a $\Bb{Z}^n$ graded 
homomorphism. From (\Rfour, Lemma 3.11) it follows that the condition 
$a_{ij} \neq a_{ik}$ for  $1\leq j,k\leq N_i$, $j\neq k$ for all 
$1\leq i\leq n$, is necessary to show that non-graded highest weight 
\lalg-module V($\cal{E}\cdot$\psila), is an irreducible 
$\underset{N- copies}{\oplus} \frk{g}$ module. Further from 
(\Rfour, Lemma 3.16), Lemma \ref{Apsi} and Lemma \ref{rank}, it follows 
that if $\lb \lambda_i \rb_{i=1}^N$ is a set of dominant integral weights, then
$A_{\text{ \psila}}$ =$\Bb{C}[t^{\Gamma}]$, for some rank n subgroup $\Gamma$ 
of $\Bb{Z}^n$. Here, $\Bb{C}[t^{\Gamma}]$ denotes the Laurent polynomial ring 
in the variables $s_i=t^{\und{m}_i}$, where $\Gamma$ is a subgroup of 
$\Bb{Z}^n$ with basis $\lb \und{m}_1,\cdots \und{m}_n \rb$.

For n=1, the map \psila\ is of the form 
\begin{equation}\begin{array}{ll} \label{num1}
\psi{\mbox{\bm{$_{(\lambda,a)}$}}}(h \otimes t_1^m) = \left(\sum_{i=1}^N \  
{\lambda}_i(h) \ a_{1i}^{m} \right)  t_1^m, \  & \forall \ h \in \frk{h},
m \in\Bb{Z}. \end{array}\end{equation}
using the following lemma, irreducible integrable modules for the 
Lie algebra $L(\frk{g}) = \frk{g} \otimes \Bb{C}[t^{\pm 1}]$ were classified
in (\CPone): 
\begin{lem}(\CPone, Lemma 4.4)\label{num} Let 
\psila : $U(\frk{h}\otimes\Bb{C}[t_1^{\pm 1}]) \rightarrow \Bb{C}[t_1^{\pm 1}]$
be a $\Bb{Z}$-graded algebra homomorphism as defined in Eqn (\ref{num1}) such 
that Image \psila = $\lar$. Then $N=rs$ for $s\in \Bb{Z}$, and there exists a 
primitive $r^{th}$ root of unity $\rho$, distinct non-zero complex numbers 
$C_{1},\cdots,C_{s}$ and a permutation $\tau$ of $\lb 1,2,\cdots,N\rb$ such 
that
\begin{eqnarray*} a_1,a_2, \cdots,a_r = C_1\rho,C_1\rho^2,\cdots,C_1\rho^r&,
\cdots ,& a_{N-r+1},a_{N-r+2}, \cdots,a_N = 
C_s\rho,C_s\rho^2,\cdots,C_s\rho^r\\
\mbox{and \qquad $\lambda_{\tau(1)}=\lambda_{\tau(2)}= \cdots =
\lambda_{\tau(r)}$} &,\cdots,& 
\lambda_{\tau(N-r+1)}= \cdots =\lambda_{\tau(N)}.
\end{eqnarray*}
\end{lem} 

\noindent The following theorem generalizes the above lemma. In the theorem, 
$pr_i:\Bb{Z}^n\rightarrow \Bb{Z}\und{e}_i$ denotes the projection map of 
$\Bb{Z}^n$ onto the $i^{th}$ coordinate. Given a rank n subgroup $\Gamma$ of
$\Bb{Z}^n$, it is easy to see that there exists positive integers 
$r_1, r_2,\cdots, r_n$ such that $pr_i(\Gamma)\cap\Gamma = r_i\Bb{Z}\und{e}_i$.

\begin{thm}\label{thm} Let \psila : $U(\frk{h}_A) \rightarrow A$ be a 
$\Bb{Z}^n$-graded algebra homomorphism as defined in Eq. (\ref{R-map}). If
Image \psila = $\Bb{C}[t^{\Gamma}]$, for a  rank n subgroup $\Gamma$ of 
$\Bb{Z}^n$ with $pr_i(\Gamma)\cap\Gamma = r_i\Bb{Z}\und{e}_i$, 
for $1\leq i\leq n$, then $r_i$ divides $N_i$ for $1 \leq$ i $\leq$ n and 
$\bf{a}$$_{N_i}$ = ($a_{i1},a_{i2}, \ldots,a_{iN_{i}}$) is a $N_i$-tuple of
distinct non-zero complex numbers  such that if ${\epsilon}_i$ is the 
${r_{_i}}^{th}$ root of unity then there exists $l_i$ ($l_i$=N$_i$/$r_{_i}$) 
distinct elements $ \lbrace c_{i1},c_{i2}, \ldots,c_{il_i} \rbrace$ in 
$\mathbb{C}^{\times}$  such that
\begin{equation}\label{a-values}
\begin{array}{lcr} a_{i1},a_{i2},\cdots, a_{ir_{_i}} &=
& {\epsilon}_i c_{i1},{\epsilon}_{i}^2 c_{i1},
\cdots,{\epsilon}_{i}^{r_{_i}} c_{i1},\\
\cdots & &\cdots\\
a_{iN_i-r_{_i}+1},a_{iN_i-r_{_i}+2},\cdots, a_{iN_i} 
&=& {\epsilon}_i c_{il_i}, {\epsilon}_{i}^2 c_{il_i},\cdots,
{\epsilon}_{i}^{r_{_i}} c_{il_i},
\end{array} \end{equation} 
and the index of $\Gamma$ in $\Bb{Z}^n$ divides N. Moreover, the irreducible
integrable \tlalg-module V(\psila) is an irreducible component of 
$\big(\otimes_{i=1}^{\ell}( V(\lambda_i)^{\otimes p})\big)\otimes A$, where
$p$=[$\Bb{Z}^n:\Gamma$].
\end{thm}
\pf  Since for each $k \in \lb 1,\cdots,n \rb$, 
$pr_k(\Gamma)\cap\Gamma = r_k\Bb{Z}e_k$ therefore,
\begin{eqnarray*}
\sum_{i_k=1}^{N_k} \left( \sum_{(i_1,\cdots,{i}_{k-1},{i}_{k+1}, \cdots,
i_n)} \lambda_{(i_1,\cdots,{i}_{k-1},i_k,{i}_{k+1}, \cdots,
i_n)} \right) \  a_{ki_{_k}}^n = 0 &
\qquad \mbox{ if \quad $n \not\equiv 0 \mod r_k$}
\end{eqnarray*}
Hence by Lemma \ref{num}, for each $i \in \lb 1,\cdots,n \rb$,
 $r_{_i}$ divides $N_i$ and there exists distinct non-zero complex numbers 
$c_{i1}, c_{i2}, \ldots, c_{il_n}$ such that Eq. (\ref{a-values}) holds, where
$ {\epsilon}_{i}$ is the $r_{i}^{th}$ root of unity and $l_i$=N$_i$/$r_{_i}$.  
But $N_i\equiv0\ {\rm{mod}}\ r_{_i}$ for each $1\leq i\leq n$, implies that
$\prod_{i=1}^n r_{_i}$ divides $\prod_{i=1}^n N_i=N$. Moreover, 
$\underset{i=1}{\overset{n}{\oplus}} r_{_i}\Bb{Z}\und{e}_i$ is a sub-lattice of
$\Gamma$, therefore p = [$\Bb{Z}^n:\Gamma$] divides 
[$\Bb{Z}^n: \underset{i=1}{\overset{n}{\oplus}} r_{_i}\Bb{Z}\und{e}_i$]=
$\prod_{i=1}^n r_{_i}$. 

The second part of the theorem is proven by inducting on n, the number of 
variables of the Laurent polynomial ring 
$\Bb{C}[t_1^{\pm 1},\cdots,t_n^{\pm 1}]$. For n=1, the proof is given by 
Lemma \ref{num}. Since  $r_i$ divides $N_i$ for $1\leq i\leq n$, therefore 
without loss of generality we may assume that $N_i=r_i$ for $1\leq i\leq n,$ 
for n$>$2. By standard theory, there exists a basis 
$\lb \und{m}^i \rb_{i=1}^n$ of $\Gamma$ such that the relation matrix 
B=($m_j^i$)$_{1\leq i,j\leq n}$ of the generators with respect to the 
standard ordered basis ($\und{e}_1,\cdots,\und{e}_n$) is lower triangular. i.e,
B is given by,
\begin{eqnarray*}B = \left(\begin{array}{cc} 
B_{n-1}& 0\\ (m_{1}^n,\cdots,m_{n-1}^n) & m^n_n \end{array}\right), 
\end{eqnarray*} where $B_{n-1}$ = ($m_j^i$)$_{1\leq i,j\leq n-1}$ is the 
relation matrix of the rank $n-1$ subgroup $\Gamma_{n-1} =
\oplus_{i=1}^{n-1}\und{m}^i\Bb{Z} $ of $\Bb{Z}^{n-1}$. Observe that for 
$1\leq i\leq n$, the element $r_i\und{e}_i \in \Gamma\cap pr_i(\Gamma)$ belongs
to the $\Bb{Z}$-linear combination of $\lb \und{m}^1,\cdots,\und{m}^i\rb$.
Hence if $\und{m}^i = (m^i_1,\cdots,m^i_i,0,\cdots,0)$, then $r_i\in 
m^i_i\Bb{Z}$, i.e., $r_i\equiv 0\mod m_i^i$ for $1\leq i\leq n$. \\
Owing to the relation matrix $B$ of $\Gamma$ we have, 
$$ {\text{\psila}}(h(\und{q}))=0\ {\rm{if}}\  
q_n\not\equiv 0\ {\rm{mod}}\ {m_n^n}.$$  Hence by Lemma \ref{num}, whenever 
$a_{n\ell_n}^{m_n^n} = a_{nk_n}^{m_n^n}$, 
\begin{equation}\begin{array}{l}\underset{(\ell_1,\cdots,\ell_{n-1})}{\sum} 
\lambda_{(\ell_1,\cdots,\ell_{n-1},\ell_n)}a_{1\ell_1}^{q_1}
\cdots a_{n-1\ell_{n-1}}^{q_{n-1}} =
\underset{(\ell_1,\cdots,\ell_{n-1})}{\sum}
\lambda_{(\ell_1,\cdots,\ell_{n-1},k_n)}a_{1\ell_1}^{q_1}\cdots 
a_{n-1\ell_{n-1}}^{q_{n-1}}. 
\end{array}\end{equation}
Since by Eqn (\ref{a-values}) for $1\leq i_n\leq N_n$, $a_{ni_n}$ is a 
$r_n^{th}$ root of a scalar and $r_n\equiv 0\mod m_n^n$, therefore given
$l_n \in \lb 1,\cdots, N_n\rb$, there exists $m_n^n$ distinct integers 
$k_n \in \lb 1,\cdots, N_n\rb$ such that $a_{n\ell_n}^{m_n^n} = 
a_{nk_n}^{m_n^n}$. Since $A_{d} = (a_{_{d i_{d}}}^{s_{d}})_{\underset{1\leq 
i_{d}\leq N_{d}}{0 \leq s_{d} \leq N_{d}-1}}$ is invertible for $1\leq d
\leq n_1$, therefore the matrix 
($a_{1i_1}^{j_1}\cdots a_{n-1i_{n-1}}^{j_{n-1}}$) = 
$\otimes_{d=1}^{n-1} A_{d}$ is invertible. Hence for a fixed (n-1)-tuple 
($\ell_1,\cdots,\ell_{n-1}$) $\in {\underset{j=1}{\overset{n-1}{\bigoplus}}} 
I_{N_j}$ there exists a subset $S_{(\ell_1,\cdots,\ell_{n-1})}$ of $I_{N_n}$ 
such that $\vert S_{(\ell_1,\cdots,\ell_{n-1})} \vert = m_n^n$ and
\[ \lambda_{(\ell_1,\cdots,\ell_{n-1},\ell_n)} = 
\lambda_{(\ell_1,\cdots,\ell_{n-1},k_n)}, \qquad \mbox{for 
$\ell_n,k_n \in S_{(\ell_1,\ell_2,\cdots,\ell_{n-1})}$.}\]
However for a fixed $\ell_n \in \lb 1,\cdots, N_n\rb$,  
\[\underset{(\ell_1,\cdots,\ell_{n-1})}{\sum} 
\lambda_{(\ell_1,\cdots,\ell_{n-1},\ell_n)}a_{1\ell_1}^{q_1}
\cdots a_{n-1\ell_{n-1}}^{q_{n-1}} =0, \quad 
\mbox{if ($q_1,\cdots,q_{n-1}$)$\notin \Gamma_{n-1}$}. \]
Thus by induction hypothesis, for a fixed $\ell_n \in I_{N_n}$, there exists 
[$\Bb{Z}^{n-1}:\Gamma_{n-1}$]= det($B_{n-1}$)
(n-1)-tuples ($\ell_1,\cdots,\ell_{n-1}$) $\in 
{\underset{j=1}{\overset{n-1}{\bigoplus}}} I_{N_j}$ for which 
$\lambda_{(\ell_1,\cdots,\ell_{n-1},l_n)}$'s are equal.
Therefore if Image \psila = $\Bb{C}[t^{\Gamma}]$ then the $\lambda_{I_j}$'s
are equal in sets of $m_n^n.[\Bb{Z}^{n-1}:\Gamma_{n-1}]$=[$\Bb{Z}^n:\Gamma$].

As a consequence, if $p$=[$\Bb{Z}^n:\Gamma$] then
 the \lalg-module $V(\cal{E}\cdot\mbox{\psila})$ is isomorphic 
to  $\underset{i=1}{\overset{M}{\otimes}} V(\lambda_i)^{\otimes p}$,
where $N=Mp$ and by Proposition \ref{irreducible comp}, 
$V(\cal{E}\cdot\mbox{\psila}) \otimes A$ is the direct sum of p irreducible 
\tlalg\ modules with V(\psila) as one of its irreducible components. 
Hence the theorem. $\hfill \blacksquare$ \endproof

Let V(\psila) be an irreducible component of 
$(\otimes_{i=1}^k V(\lambda_i))_{A}$. Then a highest weight vector of 
V(\psila) is of the form {\bm{$v$}} = $v_1\otimes\cdots\otimes v_k\otimes 
t^{\und{m}}$, where $v_i$'s are the highest weight vectors of V($\lambda_i$)'s 
for each $i$ and $t^{\und{m}} \in$ Image \psila. 
For $X\in\frk{g}\smallsetminus\frk{h}$, $\und{s}\in\Bb{Z}^n$ we have,
\begin{equation}\label{action}X(\und{s}).\mbox{\bm{$v$}}=
\sum_{i=1}^k a_{I_i}^{\und{m}} v_1\otimes\cdots
\otimes X.v_i\otimes\cdots\otimes v_k \otimes t^{\und{m}+\und{s}}. 
\end{equation} Since $ \lb X.v_1\otimes\cdots\otimes v_k\otimes t^{\und{m}}, 
\cdots, v_1\otimes\cdots\otimes X.v_i\otimes\cdots\otimes v_k\otimes 
t^{\und{m}},\cdots  \rb$ are linearly independent and  $a_{I_i}^{\und{m}}\neq 0
$ for $1 \leq i\leq N$, therefore $X(\und{s})$.{\bm{$v$}}$\neq 0$, for all 
$X \in \frk{g}\smallsetminus\frk{h}$, $\und{s}\in\Bb{Z}^n$.

\begin{prop}\label{isomorphism1}
Let $\phi,\psi : U(\hlb) \rightarrow A$ be two $\Bb{Z}^n$ graded 
algebra homomorphisms such that $\Lambda_{\psi} = \psi\vert_{\frk{h}+D}$
and  $\Lambda_{\phi+D} = \phi\vert_{\frk{h}}$. Then the \tlalg-modules 
V($\psi$) and V($\phi$) are isomorphic \iff  Ker $\psi$ = Ker $\phi$ and 
$\Lambda_{\psi} =\Lambda_{\phi}+\delta_{\und{m}}$ for some $\und{m}\in 
\Bb{Z}^n$ such that $t^{\und{m}}\in$ Image $\psi$.
\end{prop} \pf The proof is the same as in (\C, Theorem 3.5(iv)).
$\hfill \blacksquare$ \endproof

\begin{cor} \label{cor.isomorphism1}
If $\phi,\psi : U(\hlb) \rightarrow A$ be two $\Bb{Z}^n$ graded 
algebra homomorphisms such that the irreducible \tlalg-modules V($\psi$) and 
V($\phi$) are isomorphic then V($\cal{E}\cdot\psi$) is isomorphic to 
V($\cal{E}\cdot\phi$) as $\frk{g}$ modules.
\end{cor}
\pf The isomorphism of the \tlalg-modules $V(\psi)$ and $V(\phi)$ induces
an isomorphism of the \lalg-modules V($\cal{E}\cdot\psi$) and 
V($\cal{E}\cdot\phi$),
say, $f:$ V($\cal{E}\cdot\psi$) $\rightarrow$ V($\cal{E}\cdot\phi$).\\
By Lemma \ref{fin-d}, V($\cal{E}\cdot\psi$) and V($\cal{E}\cdot\phi$) are 
finite dimensional $\frk{g}$-modules. \\
Since $f$ is a \lalg-module isomorphism, 
\begin{equation}\label{4}
x_{\alpha}(\und{0}) \ f(w) =  f(x_{\alpha}(\und{0})w).
\quad \text{for $0\neq w\in V(\psi)$ and $x_{\alpha}\in\frk{g}$.}
\end{equation} 
Identifying $\frk{g}\otimes 1$ with $\frk{g}$ via the map
$x_{\alpha}(\und{0})\mapsto x_{\alpha}$, we thus see that the map $f$ induces a
$\frk{g}$-module isomorphism between V($\cal{E}\cdot\psi$) and 
V($\cal{E}\cdot\phi$). $\hfill \blacksquare$ \endproof

\para\label{isomorphism}
Recall from (\R, Theorem 1) that two $\frk{g}$ modules 
$\otimes_{i=1}^k V(\lambda_i)$ and $\otimes_{j=1}^p V(\varpi_j)$ are isomorphic
if and only if $k=p$ and $V(\lambda_i) \iso V(\varpi_{\tau(j)})$ for some 
permutation $\tau$ of $I_k$. Thus given 
({\bm{$\lambda,a$}})$\in({\mathfrak{h}}^{*})^N\times(\Bb{C}^{\times})^{N}$ and
({\bm{$\xi,b$}})$\in({\mathfrak{h}}^{*})^{N'}\times(\Bb{C}^{\times})^{N'}$,
if V(\psila) and V(\psixib) are isomorphic \tlalg\ modules, then by 
Proposition \ref{isomorphism1}, Corollary \ref{cor.isomorphism1} and 
Theorem \ref{thm}, $N=N'$, {\bm$\lambda = \xi$}$_{\tau}$ for some permutation 
$\tau$ of $\lb1,\cdots,N\rb$, suitably chosen in view of Theorem \ref{thm}.

Since Kernel \psila = Kernel \psixib, there exists a $\Bb{Z}^n$-graded 
isomorphism from Image \psila\ to Image \psixib. Hence if 
Image \psila = $\Bb{C}[t^{\Gamma}]$ then Image \psixib =  $\Bb{C}[t^{\Gamma}]$
for $\Gamma$ a subgroup of $\Bb{Z}^n$. Let $pr_i(\Gamma)\cap\Gamma = 
r_i\Bb{Z}\und{e}_i$, for $1\leq i\leq n$. For $h\in \frk{h}$ and $s\in\Bb{Z}_+,
$ $$\psi(h(sr_i\und{e}_i))\neq 0 \quad \text{if and only if}\quad  
\phi(h(sr_i\und{e}_i))\neq 0.$$ Hence if $f: V(\psi)\rightarrow V(\phi)$
is a \tlalg-module isomorphism and $v$ is the highest weight vector in $V(\psi)
$, then
\begin{eqnarray} 
h(sr_i\und{e}_i).f(v) = f(h(sr_i\und{e}_i).v), \qquad \qquad \\
\Rightarrow {\sum_{j=1}^{N_i^{^0}}(\sum_{}\xi_{_{I_{j}}}(h))b_{ii_j}^{r_{_i}s}}
f(v) = {\sum_{j=1}^{N_i^{^0}}(\sum_{}\lambda_{_{I_{j}}}(h))a_{ii_j}^{r_{_i}s}}
f(v). \label{iso.scaler} 
\end{eqnarray}
But {\bm$\lambda = \xi$}$_{\tau}$ for some permutation $\tau$ of 
$\lb1,\cdots,N\rb$, therefore Eqn (\ref{iso.scaler}) is satisfied if for each 
$i \in \lb1,\cdots,n \rb$, there exists $s_i\in\Bb{C}^{\times}$ such that 
($b_{i1},\cdots,b_{iN_i}$) = $s_i.$($a_{i\tau_i(1),\cdots,a_{i\tau_i(N_i)}}$).
 
 \noindent{\bf{Notation:}} Given V(\psila) $\in$ \Ifin , if 
$d_i$.\vpsila = $\varrho_i.$ \vpsila for all $1\leq i\leq n$, then we shall 
denote the corresponding highest weight irreducible module by
V(\psila, {\bm{$\varrho$}}). By convention we refer to V(\psila, {\bf{0}}) as 
V(\psila).

With the above notations, it follows from discussion {\bf{\ref{isomorphism}}} 
that:
\begin{thm} \label{isomorphism3}
Given ({\bm{$\lambda,a$}}) $\in (\frk{h}^*)^N \times(\Bb{C}^{\times})^N$ and
({\bm{$\xi,b$}}) $\in (\frk{h}^*)^{N_0} \times (\Bb{C}^{\times})^{N_0}$, where 
$N=N_1\cdots N_n$ and  $N_0=N_1^0\cdots N_n^0$,
and {\bm{$\varrho$}}=($\varrho_1,\cdots,\varrho_n$), 
{\bm{$\varsigma$}}=($\varsigma_1,\cdots,\varsigma_n$) $\in 
\Bb{C}^{n}$ , V(\psila, {\bm{$\varrho$}}) is isomorphic to
V(\psixib, {\bm{$\varsigma$}}) as a \tlalg-module if 
\begin{enumerate}
\item $N_i=N_i^0$ for each $1\leq i\leq n$ and consequently $N=N_0$ ;
\item For each $1\leq i\leq n$
there exists permutations $\tau_i$ of $\lb 1,\cdots,N_i \rb$ such that
\begin{enumerate}
\item $\lambda_{(i_1,\cdots,i_n)} = \xi_{(\tau_1(i_1),\cdots,\tau_n(i_n))}$ ;
\item ($b_{i1},\cdots,b_{iN_i}$) = 
$\frk{s}_i.$($a_{i\tau_i(1),\cdots,a_{i\tau_i(N_i)}}$), for some $\frk{s}_i \in
\Bb{C}^{\times}$ for $1\leq i\leq n$.\end{enumerate}
\item There exists $\und{m} \in \Gamma$ such that 
($\varrho_1,\cdots,\varrho_n$)+$\und{m}$ = ($\varsigma_1,\cdots,\varsigma_n$),
where $\Gamma$ is the subgroup of $\Bb{Z}^n$ such that Image \psila\ =
$\Bb{C}[t^{\Gamma}]$ = Image \psixib. 
\end{enumerate}
\end{thm}

\section{\small{ISOMORPHISM CLASSES OF IRREDUCIBLE \ttwist\ MODULES IN 
\Ifinmu}}\label{iso-ttwist}

In this section we shall give the isomorphism classes of the irreducible 
\ttwist\ modules in \Ifinmu. By Proposition \ref{restriction}, Proposition 
\ref{R-map} and Proposition \ref{irreducible comp'}, every irreducible 
\ttwist-module in \Ifinmu is a \ttwist\ submodule of 
V($\cal{E}\cdot$\psila)$_A$, for some $\Bb{Z}^n$ graded algebra homomorphism 
\psila. Since by (\Rfour), action of \lalg\ on the finite dimensional module 
V($\cal{E}\cdot$\psila) is given by the surjective Lie algebra homomorphism 
\begin{equation}\label{Phi(a)}
\begin{array}{c}  \Phi(\bf{a}):\frk{g}_A \rightarrow \oplus_{N-copies}
\frk{g}\\ X\otimes t_1^{m_1}\cdots t_n^{m_n}\mapsto(a_{1i_1}^{m_1}\cdots 
a_{ni_n}^{m_n} X)_{1\leq i_1\leq N_1;\cdots;1\leq i_n\leq N_n}, 
\end{array}\end{equation}
therefore action of \twist\ on V($\cal{E}\cdot$\psila) is given by the 
restriction map $\Phi(\bf{a}, \mu)= \Phi(\bf{a})\vert_{\twistb}$. It has been 
proved in (\B):
\begin{prop}(\B, Proposition 2.2) :\label{batra} Let V($\Phi(\bf{a})$)
be a finite dimensional irreducible \lalg\ module on which \lalg\ acts via the 
map $\Phi(\bf{a})$ as defined in Eq. (\ref{Phi(a)}). Then V($\Phi(\bf{a})$)
is irreducible as a \twist\ module if for each $i \in \lb 1,\cdots,n\rb$,
$\und{a}_i=(a_{i1},\cdots,a_{iN_i})$ is a $N_i$-tuple of  distinct non-zero
complex numbers and $a_{1i}^k \neq a_{1j}^k$ for $1\leq i,j\leq N_1$, 
or equivalently if Image $\Phi({\bf{a}})$ = Image $\Phi(\bf{a},\mu)$ = 
$\oplus_{N-copies} \frk{g}$. \end{prop}

\begin{rem}\label{twist red} Suppose $\oplus_{N-copies} \frk{g} \iso$ 
Image $\Phi({\bf{a}}) \neq$ Image $\Phi(\bf{a},\mu)$. Then 
Image $\Phi(\bf{a},\mu) \iso \oplus_{\cal{N}-copies} \frk{g}$ for some integer 
$\cal{N} < N$. An irreducible  $\oplus_{N-copies} \frk{g}$ module
V($\Phi(\bf{a})$)(= V($\cal{E}\cdot$\psila)) is completely reducible as a 
$\oplus_{\cal{N}-copies} \frk{g}$ module, for $\cal{N}\leq N$. Therefore any 
finite-dimensional irreducible \lalg\ module is completely reducible as a 
\twist-module under the restricted action. \end{rem}

\para \label{basis.extn}
Let Image $\Phi({\bf{a}})$ =Image $\Phi({\bf{a}}, \mu)$. Let the 
irreducible \ttwist-submodule of V($\Phi({\bf{a}}, \mu)$)$_A$ be denoted by  
V(\psilamu) where $\psi${\bm{\lamu}}:= 
$\psi{\mbox{\bm{$_{(\lambda,a)}$}}}_{\vert_{U(\thlb(\mu))}}:U(\thlb(\mu))
\rightarrow A$, denotes the restriction map. Clearly,
\[\mbox{\psila}\vert_{\frk{h}\otimes \Bb{C}[t_2^{\pm},\cdots,t_n^{\pm}]} =
\mbox{\psilamu}|_{\frk{h}\otimes\Bb{C}[t_2^{\pm},\cdots,t_n^{\pm}]}.\]
If $\lb \und{m}^1,\cdots,\und{m}^n \rb$ is a linearly independent set in
$\Bb{Z}^n$ such that Image \psila $= \Bb{C}[t^{\pm \und{m}_1},\cdots,t^{\pm 
\und{m}_n}]$ and Image \psila$\vert_{\frk{h}\otimes \Bb{C}[t_2^{\pm},\cdots,
t_n^{\pm}]}$ $= \Bb{C}[t^{\pm \und{m}_1},\cdots,t^{\pm \und{m}_{n-1}}]$,then 
the set $\lb \und{m}^1,\cdots,\und{m}^{n-1} \rb$ can be suitably extended to 
form a basis of $\Gamma^{\mu}$, where $\Gamma^{\mu} =\lb t^{\und{m}} \in A, 
t^{\und{m}}\in {\text{Image \psilamu}} \rb$. 

If V(\psilamu) is a non-trivial irreducible \ttwist-module, then  
V($\cal{E}\cdot$\psilamu) is a non-trivial \twist-module by 
Proposition \ref{irreducible comp'}. Hence the same proof as Lemma \ref{rank} 
shows that $\Gamma^{\mu}$ is a rank $n$ subgroup of $\Bb{Z}^n$. Let 
$\lb \und{m}^1,\cdots,\und{m}^{n-1},\und{\hatm}\rb$  be a basis of 
$\Gamma^{\mu}$ such that the relation matrix $B^{\mu}$ of $\Gamma^{\mu}$ 
with respect to the ordered base ($\und{e}_2,\cdots,\und{e}_n,\und{e}_1$) is 
lower triangular. 

\begin{thm}\label{thm2} Let order of $\mu$ be $k$. Assume that for 
$i \in \lb 1,\cdots,n\rb$, $\und{a}_i=(a_{i1},\cdots,a_{iN_i})$ is a 
$N_i$-tuple of  distinct non-zero complex numbers, 
$a_{1i}^k \neq a_{1j}^k$ for $1\leq i,j\leq N_1$ and for $1\leq j\leq N$,
$\lambda_{I_j}$'s are dominant integral weights. Suppose that 
Image \psila$\vert_{\frk{h}\otimes\Bb{C}[t_2^{\pm},\cdots,t_n^{\pm}]} = 
\Bb{C}[t^{\Gamma_{n-1}}]$ and $\hat{\und{m}} = \hatm_1\und{e}_2+ \cdots 
+\hatm_{n-1}\und{e}_{n}+\hatm_n\und{e}_1 \in \Bb{Z}^n$ is such that 
Image \psilamu = $\Bb{C}[t^{\Gamma^{\mu}}] = 
\Bb{C}[t^{\Gamma_{n-1}}, t^{\pm \hat{\und{m}}}]$.
Then:\begin{enumerate}
\item If $\mu(\lambda_{I_j}) \neq \lambda_{I_j}$ for some $j\in\lb1,\cdots,
N\rb$, V(\psilamu) is an irreducible \ttwist-submodule of
$(\otimes_{i=1}^\ell (V(\lambda_{I_i})^{\otimes p}))_A$, where 
$p$ =[$\Bb{Z}^{n-1}:\Gamma_{n-1}$].
\item If $\mu(\lambda_{I_j}) = \lambda_{I_j}$ for all $j\in \lb1,\cdots, N\rb$,
V(\psilamu) is an irreducible \ttwist-submodule of
$(\otimes_{i=1}^\ell (V(\lambda_{I_i})^{\otimes q}))_A$, where 
$q$ =[$\Bb{Z}^{n-1}:\Gamma_{n-1}$]$\hat{m}_n/k$. 
\end{enumerate}
\end{thm}
\noindent {\bf{Remark.}} With respect to the ordered basis 
($\und{e}_2,\cdots,\und{e}_n,\und{e}_1$) of $\Bb{Z}^n$ let 
$\und{\hatm} = (\hatm_1,\cdots,\hatm_n)$. If we prove that in case (1) 
$\hatm_n =1$ and in case (2) $\hatm_n =ks$ for some positive integer $s$,
then the desired result will follow by the same calculations as done in 
Theorem {\ref{thm}}.

\pf For $\und{m}\in\Bb{Z}^n$, if $\und{m}=(m_2,\cdots,m_n,m_1)$ with respect to
the ordered basis $\lb \und{e}_2,\cdots,\und{e}_n,\und{e}_1 \rb$ of $\Bb{Z}^n$,
then,
\begin{eqnarray}\label{twist3}
\sum_{i=1}^N \lambda_{I_i} a_{1i_1}^{m_1}a_{2i_2}^{m_2}\cdots a_{ni_n}^{m_n}
= 0, & \mbox{if $m_1\not\equiv 0\ {\rm{mod}}\ \hatm_n$ }. 
\end{eqnarray}
Since $\und{a}_i=(a_{i1},\cdots,a_{iN_i})$ are non-zero distinct complex 
numbers for each $2\leq i\leq n$, the matrix 
M =($a_{2i_2}^{m_2}\cdots a_{ni_n}^{m_n}$)$_{\underset{1\leq i_2\leq N_2;
\cdots;1\leq i_n\leq N_n}{0\leq m_2\leq N_2-1;\cdots;0\leq m_n\leq N_n-1}}$ 
being the tensor product of the invertible matrices 
$M_k= (a_{ki_k}^{m_k})_{\underset{1\leq i_k\leq N_k}{0\leq m_k\leq N_k-1}}$
for $2\leq k\leq n$, is invertible (see \B). Hence by Eq. (\ref{twist3}) 
we get, \begin{equation}\label{twist4}
\begin{array}{clcr}
\sum_{i_1=1}^{N_1} \lambda_{_{(i_1,1,\cdots,1)}}a_{1i_1}^{m_1}=0,\\
\sum_{i_1=1}^{N_1} \lambda_{_{(i_1,2,\cdots,1)}}a_{1i_1}^{m_1}=0,\\ \ldots\\
\sum_{i_1=1}^{N_1} \lambda_{_{(i_1,N_2,\cdots,N_n)}}a_{1i_1}^{m_1}=0.\\
\end{array} \quad \mbox{if} \ m_1 \not\equiv 0\ {\rm{mod}}\ \hatm_n. 
\end{equation}
Since $a_{1i}^k \neq a_{1j}^k$ for $1\leq i,j\leq N_1$, therefore 
applying (\CPtwo, Proposition 4.2), in each of the above equations we conclude 
that:

(1). $\hatm_n =1$, if $\mu({\lambda_{I_j}}) \neq \lambda_{I_j}$ for some
for some $j\in \lb1,\cdots, N\rb$. Then with notations as above,
[$\Bb{Z}^n:\Gamma^{\mu}$] = [$\Bb{Z}^{n-1}:\Gamma_{n-1}$]$\hatm_n$. Hence by 
above remark V(\psilamu) is an irreducible \ttwist-submodule of
$(\otimes_{i=1}^\ell (V(\lambda_{I_i})^{\otimes p}))_A$, where 
$p$ =[$\Bb{Z}^{n-1}:\Gamma_{n-1}$].

(2). $\hatm_1 =ks$, for some $s\geq 1$ if 
$\mu({\lambda_{I_j}}) = \lambda_{I_j}$ for all $j\in \lb1,\cdots, N\rb$.
Moreover from Eq. (\ref{twist4}) it follows that for a fixed $(n-1)-$tuple of 
integers $\hat{I}$=($i_2,i_3,\cdots,i_n$) where $1\leq i_k\leq N_k$ for 
$2\leq k\leq n$, $\lambda_{(i_1,\hat{I})}$ are equal in groups of $s$ and the 
corresponding $a_{1i_1}^k$ are proportional to the $s^{th}$ roots of unity.
Now from Proposition \ref{irreducible comp'} and a similar set of calculations 
as done in Theorem \ref{thm}, it follows that V(\psilamu) is an irreducible 
\ttwist-submodule of 
$(\otimes_{i=1}^\ell (V(\lambda_{I_i})^{\otimes q}))_A$, where 
$q$ =[$\Bb{Z}^{n-1}:\Gamma_{n-1}$]$s$,  for $s =\hat{m}_n/k$. 
$\hfill\blacksquare$\endpf

\para \label{types}
An irreducible \ttwist-module V(\psilamu) shall be referred to as:\\
$\bullet$ highest weight \ttwist-modules of first type if 
$\mu({\lambda_{I_j}}) \neq \lambda_{I_j}$ for some
for some $j\in \lb1,\cdots, N\rb$.
$\bullet$ highest weight \ttwist-modules of second type if 
$\mu({\lambda_{I_j}}) = \lambda_{I_j}$ for all $j\in \lb1,\cdots, N\rb$.

\begin{prop} \label{Twist1} Let $\Phi({\bf{a}})$, $\Phi({\bf{a}},\mu)$ and 
\psila\ be maps as  defined above. \begin{enumerate}
\item If Image $\Phi({\bf{a}})$ = Image $\Phi({\bf{a}},\mu)$ then V(\psila) is 
completely reducible as a \ttwist\ module.
\item If Image $\Phi({\bm{a}}) \neq$  Image $\Phi({\bm{a}},\mu)$ then V(\psila)
is completely reducible as a \ttwist\ module \iff\ Image {\psila} = A. 
\end{enumerate}\end{prop}
\pf  To prove (1), observe that by Proposition \ref{restriction} and 
Proposition \ref{batra},
V($\cal{E}\cdot$\psila) = V($\cal{E}\cdot$\psilamu) under the given conditions.
Thus from  Proposition \ref{irreducible comp} and 
Proposition \ref{irreducible comp'} it follows that, 
\[\mbox{V(\psila)} = \underset{\und{s}\in G^{\mu}\cap\Gamma}{\oplus} 
{\mbox{U(\ttwist).{\bm{$\omega$}}}}(\und{s}),\]
where {\bm{$\omega$}} is the highest weight vector of the irreducible \twist\ 
module V($\cal{E}\cdot$\psilamu) and given Image \psilamu = 
$\Bb{C}[t^{\Gamma^{\mu}}]$, $G^{\mu}$ is the set of coset representatives of 
$\Gamma^{\mu}$ in $\Bb{Z}^n$.

Proof of (2). Let Image $\Phi({\bm{a}}) = \underset{N-copies}{\oplus}\frk{g}$
and Image $\Phi({\bm{a}},\mu)$ = $\underset{\cal{N}-copies}{\oplus} \frk{g}$, 
with $N\neq \cal{N}$. Then by Remark \ref{twist red},
\[V(\cal{E}\cdot{\mbox{\psila}}) = \underset{i=1}{\overset{\ell}{\oplus}} 
V(\cal{E}\cdot{\mbox{\psixibmui}}), \]where 
{\bm{$\xi(i)$}} = ($\xi_{i_1},\xi_{i_2},..,\xi_{i_{\cal{N}}}$)$\in 
(\frk{h}^*)^{\cal{N}} $, {\bm{$b$}} $\in (\Bb{C}^{\times})^{\cal{N}}$ and 
$\sum_{j=1}^{\cal{N}} \xi_{i_j} \leq \sum_{j=1}^N \lambda_j$ for each 
$1\leq i\leq k$. Let 
Image \psixibmui = $\Bb{C}[t^{\Lambda_i}]$, for $1\leq i\leq k$. Then
\[(V(\cal{E}\cdot{\mbox{\psila}}))_A = \underset{i=1}{\overset{k}{\oplus}} 
\left( \underset{\und{s}_j \in G(\Lambda_i)}{\bigoplus} 
U(\ttwistb).\omega_i(\und{s}_j) \right), \]
where G($\Lambda_i$) denotes the set of coset representatives of $\Lambda_i$ in
$\Bb{Z}^n$ and $\omega_i$ is the highest weight vector of 
V($\cal{E}\cdot$\psixibmui), for $1\leq i\leq k$. 
Consider $s \in I_k$ such that 
$\sum_{j=1}^{\cal{N}} \xi_{s_j} < \sum_{j=1}^N \lambda_j$. Then
\[\omega_s = \sum {\bm{c_{s_1,..,s_N}}} 
w_{_{s_1}}\otimes... \otimes w_{_{s_N}}, \]
where weight of the each summand is $\sum_{j=1}^{\cal{N}} \xi_{s_j}$. However 
there exists a summand $w_{_{s_1}}\otimes... \otimes w_{_{s_N}}$ of $\omega_s$ 
such that $ w_{_{s_{\ell}}} \in \lb v_1, \cdots,v_{_N}\rb$ 
for some $1\leq \ell\leq \cal{N}$ but $\lb w_{_{s_1}},..,w_{_{s_N}}\rb \neq 
\lb v_1, \cdots,v_{_N}\rb$, 
where $v_i$'s are the highest weight vectors of $V(\lambda_i)$'s. In that case
it follows from  the proof of (\CPone, Proposition 1.2), that if 
$v=v_1\otimes\cdots\otimes v_N$ and Image \psila $\neq A$  then $\omega_s$ need
not be contained in U(\lalg)$v$ though U(\lalg)$v$ $\cap$ 
V($\cal{E}\cdot$\psixibmus)$_A\neq \emptyset$. Hence the proposition. 
$\hfill \blacksquare$ \endproof

\begin{notn} Given $\lambda \in \frk{h}^*$, let $\lambda^i$ =
$\lambda\vert_{\frk{h}_i}$ for $i=0,1$ if $k=2$ (resp. $i=0,1,2$ if 
$k=3$).\end{notn}

\para 
Given ({\bm{$\lambda,a$}})$\in(\frk{h}^*)^{N}\times(\Bb{C}^{\times})^N$ and 
({\bm{$\xi,b$}})$\in(\frk{h}^*)^{N'}\times(\Bb{C}^{\times})^{N'}$, such that
Image $\Phi({\bf{a}})$ = Image $\Phi({\bf{a}},\mu)$ and Image $\Phi({\bf{b}})$ 
= Image $\Phi({\bf{b}},\mu)$, the  same proof as Proposition \ref{isomorphism1}
shows that, V(\psilamu) and V(\psixibmu) are isomorphic \ttwist\ modules \iff
\begin{verse}
$\bullet$\ Ker \psilamu = Ker \psixibmu\ ;\\
$\bullet$\ \psilamu$\vert_{(\frk{h}_0+D)}$ = 
\psixibmu$\vert_{(\frk{h}_0+D)}+\delta_{\und{m}}$ for some $\und{m}\in 
\Bb{Z}^n$ such that $t^{\und{m}}\in$ Image \psilamu. 
\end{verse}
As a consequence of Proposition \ref{irreducible comp'} and Lemma \ref{twist1}
it then follows from a similar proof as Corollary \ref{cor.isomorphism1} that 
if V(\psilamu) and V(\psixibmu) are isomorphic as \ttwist\ modules then
V($\cal{E}\cdot$\psilamu) is isomorphic to V($\cal{E}\cdot$\psixibmu) as 
$\frk{g}$ modules.    Consequently,  by (\R, Theorem 1), $N=N'$ 
and {\bm{$\lambda$}}={\bm{$\xi$}}$_{\tau}$, for some $\tau$ in the permutation 
group of $I_N$ suitably chosen in view of Theorem \ref{thm2} and 
Image \psilamu\ = $\Bb{C}[t^{\Gamma^{\mu}}]$ = Image \psixibmu\ 
for some rank n subgroup $\Gamma^{\mu}$ of $\Bb{Z}^n$.

Let $A_0 = \Bb{C}[t_1^{\pm 2},t_2^{\pm},\cdots,t_n^{\pm}]$. 
Since $\Bb{C}[t_2^{\pm},\cdots,t_n^{\pm}] \subset A_0$, if 
Image \psilamu$\vert_{(\frk{h}\otimes\Bb{C}[t_2^{\pm},\cdots,t_n^{\pm}])(\mu)} 
= \Bb{C}[t^{\Gamma^{\mu}_{n-1}}]$
then $t^{\und{m}}\in A_0$ for all $\und{m} \in \Gamma^{\mu}_{n-1}$.

With notation as in Proposition \ref{isomorphism3} we prove the following 
theorem:
\begin{thm}\label{tiso} Let V(\psilamu, {\bm{$\varrho$}}) and 
V(\psixibmu, {\bm{$\varsigma$}}) be two  irreducible \ttwist-modules in \Ifinmu
such that Image $\Phi({\bf{a}})$ =  Image $\Phi({\bf{a}},\mu)$ and Image 
$\Phi({\bf{b}})$ =  Image $\Phi({\bf{b}},\mu)$. Then they are isomorphic if:
\begin{enumerate}
\item Both are of the same type. (see \ref{types})
\item If both are of the first type then for each $1\leq i\leq n$ there exists 
permutations $\tau_i$ of $\lb 1,\cdots,N_i \rb$ such that
\begin{enumerate}
\item ($b_{i1},\cdots,b_{iN_i}$) = $\wp_i.$($a_{i\tau_i(1),\cdots,
a_{i\tau_i(N_i)}}$) for $2\leq i\leq n$ and  \\
($b_{11}^k,\cdots,b_{1N_1}^k$)= $\wp_1.$($a_{i\tau_i(1)}^k,\cdots,
a_{i\tau_i(N_i)}^k$), i.e., for each $1\leq i\leq N_1$, $b_{1i}= 
\varepsilon \wp_1 a_{1\tau_1(i)}$, where  $\varepsilon$ is a $k^{th}$ root of 
unity and $\wp_i \in  \Bb{C}^{\times}$ for $1\leq i\leq n$.
\item $\xi_{(i_1,\cdots,i_n)}^0 = \lambda_{(\tau_1(i_1),\cdots,\tau_n(i_n))}^0$
and $\xi_{(i_1,\cdots,i_n)}^i = \varepsilon^{-i} \lambda_{(\tau_1(i_1),
\cdots, \tau_n(i_n))}^i $ (for i=1 when k=2 and i=1,2 when k=3)  whenever 
$b_{1i_1} = \varepsilon\wp_1 a_{1\tau_1(i_1)}$,
where $\varepsilon$ is a $k^{th}$ root of unity .
\end{enumerate}
\item If both are of the second type then for all $1\leq i\leq n$, $N_i=N_i^0$ 
and there exists  permutations  $\tau_i$ of $\lb 1,\cdots,N_i \rb$ such that 
\begin{enumerate}
\item $\xi_{(i_1,\cdots,i_n)}= \xi_{(i_1,\cdots,i_n)}^0 = 
\lambda_{(\tau_1(i_1),\cdots,\tau_n(i_n))}^0 = 
\lambda_{(\tau_1(i_1),\cdots,\tau_n(i_n))}.$
\item ($b_{i1},\cdots,b_{iN_i}$) = $\wp_i.$($a_{i\tau_i(1)},\cdots,
a_{i\tau_i(N_i)}$) for $2\leq i\leq n$ and  \\
($b_{11}^k,\cdots,b_{1N_1}^k$)= $\wp_1.$($a_{i\tau_i(1)}^k,\cdots,
a_{i\tau_i(N_i)}^k$), i.e., for each $1\leq i\leq N_1$, $b_{1i}= 
\varepsilon \wp_1 a_{1\tau_1(i)}$, where  $\varepsilon$ is a $k^{th}$ root of 
unity and $\wp_i \in  \Bb{C}^{\times}$ for $1\leq i\leq n$. \end{enumerate}
\item There exists $\und{m} \in \Gamma^{\mu}$  such that 
($\varrho_1,\cdots,\varrho_n$)+$\und{m}$ = ($\varsigma_1,\cdots,\varsigma_n$).
where $\Gamma^{\mu}$ is a rank n subgroup of $\Bb{Z}^n$ such that 
Image \psilamu =$\Bb{C}[t^{\Gamma^{\mu}}]$= Image \psixibmu.
\end{enumerate} \end{thm}
\pf (1). Let V(\psilamu) and V(\psixibmu) be irreducible \ttwist\ modules of 
first and second type respectively.
We consider the case when k=2. The case k=3, can be proved similarly.  

For k=2, note that by \ref{basis.extn}, there exists $h \in \frk{h}_1$ such 
that \[ \mbox{\psilamu}(h\otimes t_1t^{\und{m}}) \neq 0, \qquad \mbox{for some 
$t^{\und{m}} \in A_0$.}\] 
But \psixibmu$\vert_{\frk{h}\otimes A_1}=0$. Hence Ker \psilamu $\neq$ 
Ker \psixibmu. Hence (1). 

\noindent (2). Let V(\psilamu) and  V(\psixibmu) be isomorphic irreducible 
\ttwist\ modules of first type. Since \twist = 
$\oplus_{i=0}^{k-1} (\frk{g}_i\otimes t^i_1A_0)$, for k=2,3,  
where $\frk{g}_0$ is a simple finite dimensional Lie algebra, the 
\ttwist-modules V(\psilamu) and V(\psixibmu) are also isomorphic as
modules for the multi-loop Lie algebra $\frk{g}_0 \otimes A_0$. Consequently, 
by Theorem \ref{isomorphism3}, 
$N_i=N_i'$ for each $1\leq i\leq n$ and 
there exists  permutations $\tau_i$ of $\lb 1,\cdots,N_i \rb$ 
such that
$$ \begin{array}{ll}
 \xi_{(i_1,\cdots,i_n)}^0 = \lambda_{(\tau_1(i_1),\cdots,\tau_n(i_n))}^0,&\\ 
(b_{i1},\cdots,b_{iN_i}) = \wp_i.(a_{i\tau_i(1),\cdots,a_{i\tau_i(N_i)}}), 
\qquad &\text{for $2\leq i\leq n$},  \\
(b_{11}^2,\cdots,b_{1N_1}^2) = 
\wp_1.(a_{i\tau_i(1)}^2,\cdots,a_{i\tau_i(N_i)}^2), \qquad 
&\text{where $\wp_i \in  \Bb{C}^{\times}$ for $1\leq i\leq n$}. 
\end{array}$$ 
Since $b_{1i}^2 = \wp_1 a_{1\tau_1(i)}^2$, therefore 
$b_{1i}= \pm \sqrt{\wp_1} a_{1\tau_1(i)}$ for each $1\leq i\leq N_1$.
But we have seen that {\bm{$\lambda$}}={\bm{$\xi$}}$_{\tau}$, for some 
$\tau$ in the permutation group of $I_N$, hence there exists permutations 
$\tau_i$ of $\lb 1,\cdots,N_i \rb$ such that $\xi_{(i_1,\cdots,i_n)} = 
\lambda_{(\tau_1(i_1),\cdots,\tau_n(i_n))}$. However, for 
$h \in \frk{h}_j$, $j\neq 0$, and $t^{\und{m}} \in A_0$, 
\[\mbox{\psilamu}(h\otimes t_1^jt^{\und{m}}) = 
\left(\sum_{i=1}^n \lambda_{I_i}^j 
a_{1i_1}^{km+j}a_{2i_2}^{m_2}\cdots a_{ni_n}^{m_n}\right) t_1^jt^{\und{m}}.\] 
Hence for k=2,\ $\xi_{(i_1,\cdots,i_n)} = \lambda_{(\tau_1(i_1),\cdots,
\tau_n(i_n))}$, \ 
$\xi_{(i_1,\cdots,i_n)}^0 = \lambda_{(\tau_1(i_1),\cdots,\tau_n(i_n))}^0$, 
where $\tau_i$ is a permutation of $\lb 1,\cdots,N_i \rb$ for $1\leq i\leq n$,
and \ $b_{1i}= \pm \sqrt{\wp_1} a_{1\tau_1(i)}$ for each $1\leq i\leq N_1$ 
imply:
\begin{eqnarray*}
 \xi_{(i_1,\cdots,i_n)}^1 = \pm \lambda_{(\tau_1(i_1),\cdots,\tau_n(i_n))}^1, 
\qquad {\rm{whenever}} \ \ 
b_{1i_1}= \pm \wp_1 a_{1\tau_1(i_1)}. \end{eqnarray*}Hence (2) holds for $k=2$.
The case k=3, can be proved similarly. \\
\noindent(3). Let V(\psilamu) and V(\psixibmu) be two isomorphic irreducible 
\ttwist modules of the second type. Suppose that Image \psilamu = 
$\Bb{C}[t^{\Gamma^{\mu}}]$ = Image \psixibmu. Then from 4.3 it is clear 
that for all  $\und{m} \in \Gamma^{\mu},\ t^{\und{m}} \in A_0.$ Hence by the 
same analysis as above (3) can be proved. \\
\noindent(4) is clear.  $\hfill \blacksquare$ \endpf

\subsection*{\large{ACKNOWLEDGMENTS}}
The authors are grateful to Prof. S.Eswara Rao for suggesting the problem and 
pointing out errors in an earlier draft of the manuscript. The first author 
would also like to thank Prof Amritanshu Prasad, Dr. Anupam Kumar Singh
and Prof. Maneesh Thakur for helpful discussions during this work. She would 
also like to thank Dr. Purusottam Rath for bringing reference (\R) to her 
notice.

\subsection*{\large{REFERENCES}}

\noindent Allison, B., Berman, S., Faulkner, J., Pianzola, A. Realization of 
  graded simple algebras as loop algebras, arXiv:math/0511723v2 [math.RA],to 
  appear in Forum Mathematicum.

\noindent Batra, P. (2004). Representations of twisted multi-loop Lie algebras.
    J. Algebra 272:404-416.

\noindent Chari, V. (1986). Integrable representations of affine Lie-algebras. 
      Invent. Math. 85:317-335.

\noindent Chari, V., Pressley, A. N. (1986). New unitary representations of 
loop groups.  Math. Ann. 275:87-104.

\noindent Chari, V., Pressley, A. N. (1988) Integrable representations of 
twisted affine Lie algebras. J. Algebra 113:438-464.

\noindent Humphreys, J. E. (1972). Introduction to Lie-algebras and
representations theory.  Springer, Berlin, Hidelberg, New York.

\noindent Kac, V. G. (1990). Infinite dimensional Lie algebras, 3rd ed. 
Cambridge University Press.  

\noindent Rajan, C. S. (2004). Unique decomposition of tensor products of 
irreducible representations of simple algebraic groups. Ann. of Math. 
160:683-704.

\noindent Rao, S. Eswara. (1993). On Represetations of loop algebra. 
Comm. Algebra 21(6):2131-2153. 

\noindent Rao, S. Eswara. (1995). Iterated loop modules and a filtration for 
the vertex representation of toroidal Lie algebras. 
Pacific J. Math 171(2):511-528.

\noindent Rao, S. Eswara. (2001). Classification of irreducible integrable 
modules for multi-loop Lie algebras with finite dimensional weight spaces. 
J. Algebra 246:215-225.

\noindent Rao, S. Eswara. (2004). Classification of irreducible integrable 
modules for toroidal Lie algebras with finite dimensional weight spaces. 
J. Algebra 277:318-348.

\noindent Yoon, Y. (2002). On polynomial representations of current algebras.
     J. Algebra 252(2):376-393.

\end{document}